\documentclass[11pt]{amsart} 

\usepackage{marvosym}
\usepackage{amsmath}
\usepackage{amssymb}
\usepackage{amsthm}
\usepackage{color}
\usepackage{graphicx}
\usepackage{pst-all}
\usepackage{pb-diagram}
\usepackage{stmaryrd }
\usepackage{url}
\usepackage{bm}

\usepackage[utf8]{inputenc}
\usepackage{enumitem}
\usepackage{enumitem}
\usepackage[T1]{fontenc}
\usepackage{verbatim}
\usepackage{tikz}
\usepackage{pgfplots}
\usepackage{caption} 
\usepackage{hyperref}
\hypersetup{
   colorlinks,
   citecolor=red,
   filecolor=black,
   linkcolor=blue,
   urlcolor=black
 }%

\newcommand{\bigzero}{\mbox{\normalfont\Large\bfseries 0}}
\newcommand{\rvline}{\hspace*{-\arraycolsep}\vline\hspace*{-\arraycolsep}}

\newtheorem{prop}{Proposition}[section]
\newtheorem{lema}[prop]{Lemma}
\newtheorem{teor}[prop]{Theorem}
\newtheorem{fact}[prop]{Fact}
\newtheorem{obs}[prop]{Observation}
\newtheorem{preg}[prop]{Question}
\newtheorem{corol}[prop]{Corollary}
\newtheorem{conj}[prop]{Conjecture}

\theoremstyle{definition}
\newtheorem{ejem}[prop]{Example}
\newtheorem{rmk}[prop]{Remark}
\newtheorem{const}[prop]{Construction}
\newtheorem{definicion}[prop]{Definition}
\newtheorem{notation}[prop]{Notation}

\def\Ind#1#2{#1\setbox0=\hbox{$#1x$}\kern\wd0\hbox to 0pt{\hss$#1\mid$\hss}     
\lower.9\ht0\hbox to 0pt{\hss$#1\smile$\hss}\kern\wd0}
\def\Notind#1#2{#1\setbox0=\hbox{$#1x$}\kern\wd0\hbox to 0pt{\mathchardef
\nn=12854\hss$#1\nn$\kern1.4\wd0\hss}\hbox to
0pt{\hss$#1\mid$\hss}\lower.9\ht0 \hbox to
0pt{\hss$#1\smile$\hss}\kern\wd0}

\def\ind{\mathop{\mathpalette\Ind{}}}                              % FOR INDEPENDENCE SYMBOL
\newcommand{\bp}{\begin{prop}}
\newcommand{\ep}{\end{prop}}
\newcommand{\bd}{\begin{definicion}}
\newcommand{\ed}{\end{definicion}}
\newcommand{\bej}{\begin{ejem}}
\newcommand{\eej}{\end{ejem}}
\newcommand{\bl}{\begin{lema}}
\newcommand{\el}{\end{lema}}
\newcommand{\bh}{\begin{fact}}
\newcommand{\eh}{\end{fact}}
\newcommand{\bpreg}{\begin{preg}}
\newcommand{\epreg}{\end{preg}}
\newcommand{\bo}{\begin{obs}}
\newcommand{\eo}{\end{obs}}
\newcommand{\bcon}{\begin{conj}}
\newcommand{\econ}{\end{conj}}
\newcommand{\brmk}{\begin{rmk}}
\newcommand{\ermk}{\end{rmk}}
\newcommand{\bc}{\begin{corol}}
\newcommand{\ec}{\end{corol}}
\newcommand{\bconst}{\begin{const}}
\newcommand{\econst}{\end{const}}
\newcommand{\bdem}{\begin{proof}}
\newcommand{\edem}{\end{proof}}
\newcommand{\benum}{\begin{enumerate}}
\newcommand{\eenum}{\end{enumerate}}
\newcommand{\bitem}{\begin{itemize}}
\newcommand{\eitem}{\end{itemize}}
\newcommand{\be}{\begin{ejem}}
\newcommand{\ee}{\end{ejem}}
\newcommand{\bt}{\begin{teor}}
\newcommand{\et}{\end{teor}}

\newcommand{\Los}{\L o\'s}

\newcommand{\G}{\mathcal{G}}

\newcommand{\Z}{\mathbb{Z}}

\newcommand{\Nat}{\mathbb{N}}
\newcommand{\Q}{\mathbb{Q}}

\newcommand{\ov}[1]{\overline{#1}}
\newcommand{\Le}{\mathcal{L}}
\newcommand{\U}{\mathcal{U}}

\def\bigcupdot{\mathop{\hbox to 0pt{\hskip.43em$\dot{\ }$\hss}\bigcup}}

\definecolor{electricindigo}{rgb}{0.44, 0.0, 1.0}
\definecolor{ao(english)}{rgb}{0.0, 0.5, 0.0}
\definecolor{dunkelmagenta}{rgb}{.7, 0, 0}
\definecolor{orange}{rgb}{0.6,0.2,0}

\title{Pseudofinite abelian groups and asymptotic classes of finite abelian groups}
\author{Carlos C\'ordoba-Caycedo, Dar\'io Garc\'ia}
\date{}
\address{Carlos C\'ordoba Caycedo \\
Department of Mathematics $\&$ Statistics, McMaster University, Canada}
\email{cordobac@mcmaster.ca}
\address{Dar\'io Garc\'ia  \\
Departamento de Matem\'aticas \\
Universidad de los Andes \\
Carrera 1 No. 18A-10, Bogot\'a D.C.}
\email{da.garcia268@uniandes.edu.co}

\begin{document}

\begin{abstract} This paper studies pseudofinite abelian groups and counting properties of definable sets in finite abelian groups and their ultraproducts. First, we characterize pseudofinite abelian groups in terms of Szmielew invariants and, building on results of Macintyre, extend this analysis to obtain a characterization of pseudofinite $\omega$-stable abelian groups. We then give an alternative proof that the class of finite abelian groups forms a multidimensional exact class, explicitly describing the possible cardinalities of definable sets in classes of finite abelian groups.\\

We also investigate $N$-dimensional asymptotic classes of finite abelian groups and provide a series of examples showing that, despite the exactness result, several natural classes of finite abelian groups fail to be $N$-dimensional asymptotic classes for every $N\geq 1$. Finally, in an appendix, we characterize the SU-rank of pseudofinite abelian groups, building on work of Rogers. In particular, we show that for every pseudofinite abelian group $A$, either $\operatorname{Th}(A,+)$ is superstable of finite SU-rank or it is strictly stable.
\end{abstract}

\maketitle
\tableofcontents

\section{Introduction}

A \emph{pseudofinite structure} is an $\Le$-structure $M$ such that every
$\Le$-sentence that holds in $M$ has a finite model, or equivalently, an
$\Le$-structure that is elementarily equivalent to an ultraproduct of finite
$\Le$-structures. Among abelian groups, the simplest examples are the infinite vector
spaces over $\mathbb{F}_p$ and the group $(\mathbb{R},+)$, which is
isomorphic to any non-principal ultraproduct of the groups in the class
$\{(\mathbb{Z}/p\mathbb{Z},+) : p \text{ prime}\}$. A less obvious example,
which we will see repeatedly, is the group
$\mathbb{Z}(p^\infty)\oplus\mathbb{Z}_{(p)}$, the direct sum of the \emph{Prüfer $p$-group} with the localization of $\mathbb{Z}$ in the ideal $(p)$, which is elementarily
equivalent to a non-principal ultraproduct of the finite cyclic groups
$\mathbb{Z}/p^k\mathbb{Z}$ (see Lemma \ref{SumTwoNoPseudofiniteGroups}).\\

This paper studies which abelian groups are pseudofinite, and which measurability conditions
are satisfied by ultraproducts of finite abelian groups. The model theory of abelian groups is unusually well understood. By the
work of W. Szmielew (see \cite[Appendix A.2]{Hodges-ModelTheory}), the complete theory of an abelian group $A$ is determined by a countable
family of numerical invariants, the \emph{Szmielew invariants}
$U(p,n;A)$, $D(p,n;A)$, $\mathrm{Tf}(p,n;A)$ and $\mathrm{Exp}(p,n;A)$, which
record the sizes of quotients of definable subgroups. A. Macintyre in \cite{Macintyre1971} used this machinery to characterize the $\omega$-stable
abelian groups, and P. Rogers in \cite{PhD-PatriciaRogers} characterized the superstable
ones.\\

On the other hand, a description of the \emph{pseudofinite} abelian
groups in terms of these invariants is, to our knowledge, not available in
a convenient reference. Perhaps the best approximation is Theorem A.2.14 of \cite{Hodges-ModelTheory}, which itself refers to a paper of Basarab from 1971 and another of Ershov from 1963, but despite our best efforts we could find neither of these papers. In
Section \ref{sect:Characterizationpsfgroups} we fill this gap by showing that an abelian group $G$ is
pseudofinite if and only if $D(p,n;G)=\mathrm{Tf}(p,n;G)$ for every prime $p$ and every $n$, an equality that holds in every finite abelian group, and we
translate this into an explicit condition on the parameters of a Szmielew
group (Theorem \ref{Theorem:ClassificationPsfAbelianGroups}). Along the way we show that the
class of pseudofinite abelian groups is closed under countable direct sums, and based on  Macintyre's theorem, we
then characterize the abelian groups that are both pseudofinite and
$\omega$-stable, and we show with explicit examples that these two properties
are independent (Theorem \ref{ClassificationWStablePsfAbGroups}).\\

One of the main features of ultraproducts of finite structures is the
analysis of the non-standard cardinalities of definable sets. In the
general case very little can be said about these cardinalities, but under strong
conditions on the class of finite structures it is possible to describe
them precisely. Macpherson and Steinhorn introduced in \cite{Macpherson-Steinhorn}  the concept of 1-dimensional asymptotic classes, Elwes generalized it in \cite{Elwes-AsymptoticClasses} to
\emph{$N$-dimensional asymptotic classes}, and more recently Anscombe, Macpherson, Steinhorn and
Wolf \cite{AMSW} introduced the concept of multidimensional asymptotic classes (m.a.c.) and \emph{multidimensional exact
classes} (m.e.c.), which are flexible enough to include classes of finite structures whose
ultraproducts have infinite rank. Some prominent examples include the class of
finite fields, the class of finite cyclic groups (both in \cite{Macpherson-Steinhorn}), classes whose ultraproducts are strongly minimal (see \cite{Pillay-stronglyminimalpsf}) or uncountably
categorical (cf. \cite{vanAbel-counting}), and the classes of regular graphs of growing
degree and girth studied in \cite{Garcia-Robles}.\\

In this paper we focus on measurability conditions for classes of finite
abelian groups. It was shown in \cite{AMSW} that the class of all finite
abelian groups is a multidimensional exact class, as a consequence of the
analogous result for finite $R$-modules, although without an explicit
description of the possible cardinalities of definable sets. In
Section \ref{sect:Characterizationpsfgroups} we give an alternative proof, based on the Smith normal form and on
the structure of cosets of definable subgroups, which shows that the cardinality of a set definable by a formula $\varphi(x;\ov{a})$ in one variable has only finitely many possibilities depending on cardinalities of the form $|A|^{1+\mu-\nu}$, $\{|A/c_\ell A|:\ell\leq\nu\}$ and $\{|A/d_\ell A|:\ell\leq\mu\}$, where the integers $\mu,\nu,c_\ell,d_\ell$
depend only on the formula $\varphi(x;\ov{y})$. (Lemma \ref{Lemma:sizes-pp-definable}). The Projection Lemma of Wolf then extends the result to formulas in
several variables.\\

Another natural notion of measurability is given in \cite{Macpherson-Steinhorn} by 
Macpherson and Steinhorn, and its asymptotic counterparts. MS-measurability in theories of modules was studied by Kestner in \cite{Kestner2014}. It is claimed in 
\cite[Remark 4.14]{Kestner2014} that the MS-measurable abelian groups are
precisely the pseudofinite ones. In Section \ref{SS:MS-measurability} we show that this is not
the case: the ultraproduct of the finite groups
$\mathbb{Z}/p\mathbb{Z}\oplus\cdots\oplus\mathbb{Z}/p^t\mathbb{Z}$ is
pseudofinite but not superstable, hence not MS-measurable, since
MS-measurable structures are supersimple of finite SU-rank. \\

Section \ref{SS:AsymptoticClasses} addresses the problem of determining which classes of finite
abelian groups are $N$-dimensional asymptotic classes. The finite cyclic
groups are a 1-dimensional asymptotic class \cite[Theorem 3.14]{Macpherson-Steinhorn}, but we show with several examples that despite the exactness result several natural
classes of finite abelian groups (for instance those closed under direct
sums, or of the form $\{\mathbb{Z}/p\mathbb{Z}^{\,m}\oplus
\mathbb{Z}/q\mathbb{Z}^{\,n}\}$ for two fixed primes) fail to be
$N$-dimensional asymptotic classes for every $N\ge1$.\\

The obstruction is always the same: ultraproducts of $N$-dimensional
asymptotic classes are supersimple of finite SU-rank \cite{Elwes-AsymptoticClasses}, and since a stable supersimple
theory is superstable, Rogers' criterion shows
that in both cases the set $\Sigma$ of pairs $(p,n)$ with
$|p^nA/p^{n+1}A|=\infty$ must be finite. The key tools are then a criterion
in terms of the width and depth of the groups in the class
(see Lemmas \ref{Lemma:wdinfty} and \ref{lem:InfiniteWidthforInfPrimes}) that produces
non-superstable ultraproducts. We close with one conjecture on the general picture as an attempt to characterize those classes of finite abelian groups corresponding to $N$-dimensional asymptotic classes. (see Conjecture \ref{Conj:NDimAsympt})\\

Finally, in the Appendix we characterize forking and SU-rank in arbitrary pseudofinite
abelian groups, building on the work of Rogers. Writing
$\Sigma(A)=\{(p,n) : |p^nA/p^{n+1}A|=\infty\}$, Rogers shows that $\mathrm{Th}(A)$
is superstable if and only if $\Sigma(A)$ is finite, and when $A$ is pseudofinite and superstable, the SU-rank of $A$ is equal to $|\Sigma(A)|$ or $|\Sigma(A)|+1$ depending on whether
$tA$ is finite, where $t$ is the least common multiple of the relevant prime
powers. In particular, for every pseudofinite abelian group $A$, either $\mathrm{Th}(A)$
is superstable of finite SU-rank or it is strictly stable. These results
are used in several examples of Section \ref{SS:AsymptoticClasses}.

\subsection*{Acknowledgements} The first author is deeply grateful for the financial support of the Fields Institute in the writing of this paper.

\section{Preliminaries}\label{sec:prelim}

This section collects the background material used throughout the rest of the paper. In Subsection \ref{Subsec:BasicFacts} we summarize some elementary facts about finite abelian groups, together with the Smith Normal Form. Afterwards, in Subsection \ref{Subsec:ModelTheoryAG} we turn to the model-theoretic side: we review Prüfer's Lemma, which reduces the study of the complete theory of an abelian group to its positive-primitive formulas, and then Szmielew's theorem, which shows that every abelian group is uniquely determined by a family of numerical invariants, called the Szmielew invariants.

\subsection{Basic facts about finite abelian groups}
\label{Subsec:BasicFacts}We use this subsection to give a list of easy results about finite abelian groups that will be used throughout the paper. Most of the results are fairly straightforward, but we point out in particular the Notation \ref{Notation-Groups} which may be a bit different from the standard notation on abelian groups.

\bl \label{Lemma:intersection-cosets} Let $k,\ell$ be two integers such that $\gcd(k,\ell)=1$. Then every coset in $A/kA$ intersects every coset of $A/\ell A$.
\el
\bdem To see this, let us assume that $a,b\in A$ are two arbitrary elements of $A$. It suffices to find an element $c\in A$ such that $c-a\in kA$ and $c-b\in \ell A$, as this implies that the coset $a+kA$ intersects the coset $b+\ell A$.\\

Since $\gcd(k,\ell)=1$, there are integers $r,s$ such that $r\cdot k+s\cdot \ell=1$. Consider the element $c=\ell s\cdot a+rk\cdot b$, in order to obtain:
\begin{align*}
c-a&=\ell s\cdot a+kr\cdot b-a=(1-rk)\cdot a+kr\cdot b-a=kr\cdot (-a+b)\in kA,\\
c-b&=\ell s\cdot a+kr\cdot b-b=\ell s\cdot a+(1-\ell s)\cdot b-b=\ell s\cdot (a-b)\in \ell A. \qedhere
\end{align*}
\edem

\bl \label{Lemma:intersections-lcm} Let $A$ be an abelian group, and $k,\ell$ positive integers with $r=\operatorname{lcm}(k,\ell)$. Then, we have  $(kA)\cap (\ell A)=rA$.
\el
\bdem We can write $r=k\cdot m_1=\ell\cdot m_2$ for some $m_1,m_2\in\mathbb{N}$, and so we have $rA\subseteq (k\cdot m_1)\cdot A=k(m_1A)\subseteq kA$ and $rA\subseteq (\ell\cdot m_2)A=\ell(m_2A)\subseteq \ell A$, and so $rA\subseteq kA\cap \ell A$.\\

Conversely, suppose that $a\in kA\cap \ell A$, and so there are $b_1,b_2\in A$ such that $a=k\cdot b_1=\ell\cdot b_2$. Let us write $d=\gcd(k,\ell)$, and recall that $k\cdot \ell=\dfrac{\operatorname{lcm}(k,\ell)}{\gcd(k,\ell)}=\dfrac{r}{d}$, and recall that there are integers $\alpha,\beta$ such that $d=k\cdot \alpha+\ell\cdot \beta$. Hence, we can consider the element $c=\alpha\cdot b_2+\beta\cdot b_1$, and so we have:
\begin{align*}
r\cdot c&=r\cdot(\alpha\cdot b_2+\beta\cdot b_1)=\dfrac{k\cdot \ell}{d}(\alpha\cdot b_2+\beta\cdot b_1)\\
&=\dfrac{k\cdot \alpha}{d}(\ell\cdot b_2)+\dfrac{\ell\cdot\beta}{d}(k\cdot b_1)=\dfrac{k\cdot \alpha}{d}\cdot a+\dfrac{\ell\cdot\beta}{d}\cdot a\\
&=\dfrac{k\cdot \alpha+\ell\cdot \beta}{d}\cdot a=\dfrac{d}{d}\cdot a=a. 
\end{align*}
Therefore, $c\in rA$, and so we conclude that $kA\cap \ell A\subseteq rA$.
\edem

\bl \label{Lemma-IntersectionCosets} Let $A$ be an abelian group. If $G,H$ are subgroups of $A$ and $c_1,c_2\in A$, then the intersection of cosets $(c_1+G)\cap (c_2+H)$ is either empty or a coset of the subgroup $G\cap H$.
\el
\bdem Suppose the intersection is non-empty, and let $b$ be an arbitrary element in $(c_1+G)\cap (c_2+H)$. Then, there are elements $g\in G$ and $h\in H$ such that $b=c_1+g=c_2+h$. Let us now show that $(c_1+G)\cap (c_2+H)=b+(G\cap H)$.\\

If $\alpha\in (c_1+G)\cap (c_2+H)$, there are elements $g_\alpha\in G$, $h_\alpha\in H$ such that $\alpha=c_1+g_\alpha=c_2+h_\alpha$. Hence, $\alpha-b=g_\alpha-g=h_\alpha-h$, which implies that $\alpha-b\in (G\cap H)$. Hence, $\alpha\in b+(G\cap H)$.\\

Conversely, if $\beta\in b+(G\cap H)$, then there is an element $t\in G\cap H$ such that $\beta=b+t$. Using the identities for $b$, we obtain $\beta=b+t=c_1+(g+t)$ and $\beta=c_2+(h+t)$, which shows that $\beta\in (c_1+G)\cap(c_2+H)$. 
\edem

\bc \label{Corollary:IntersectionCosets} Let $A$ be an abelian group. If $a_1+G_1,\ldots,a_k+G_k$ are finitely many cosets of subgroups of $A$, then either the intersection $\displaystyle{\bigcap_{i=1}^k (a_i+G_i)}$ is empty, or it is equal to a coset of the form $\displaystyle{b+\bigcap_{i=1}^kG_i}$. 
\ec
\bdem This follows directly from Lemma \ref{Lemma-IntersectionCosets}. Namely, in the inductive step we can take $G=G_1\cap\cdots\cap G_k$ and $H=G_{k+1}$.
\edem

\bd Let $A$ be an abelian group. For $m\in \mathbb{Z}^{\geq 0}$, we can define the function $\lambda_m:A\to A$ by $\displaystyle{\lambda_m(x):=m\cdot x=\underbrace{x\ +\ \cdots\ +\ x}_{\text{$m$-times}}.}$ Similarly, we can define  $\lambda_{(-m)}(x):=-\lambda_{m}(x)$.
\ed

\bl \label{Lemma:preimage-element} Let $A$ be an abelian group, $a\in A$ and $k\in\mathbb{Z}$. Then, the solution set of the equation $k\cdot x+a=0$ is either empty, or it has the form $\alpha_0+\operatorname{ker}(\lambda_k)$ for some (any) solution $\alpha_0$ from the equation. 
\el

\bdem Suppose that $\alpha_0$ is an arbitrary solution of the equation $k\cdot x+a=0$. We will show that for $\{\alpha\in A: k\cdot \alpha+a=0\}=\alpha_0+\operatorname{ker}(\lambda_k)$. Since $\alpha_0$ is a solution of the equation, we have $k\cdot \alpha_0+a=0$. So, if $\alpha\in (\alpha_0+\operatorname{ker}(\lambda_k))$, then $\lambda_k(\alpha-\alpha_0)=0$, which implies $k\cdot \alpha=k\cdot \alpha_0$. Hence, we have $k\cdot \alpha+a=k\cdot \alpha_0+a=0$, and so $\alpha$ is a solution of the equation $k\cdot x+a=0$.\\

Conversely, if $\alpha$ is a solution of the equation $k\cdot x+a=0$, then we have 
$k\cdot \alpha=-a=k\cdot \alpha_0$, so $k(\alpha-\alpha_0)=0$. That is, $\alpha-\alpha_0\in \operatorname{ker}(\lambda_k)$.
\edem

\bl \label{Lemma:intersection-kernels} Let $A$ be an abelian group and $m_1,\ldots,m_k\geq 1$. Then, for $g=\gcd(m_1,\ldots,m_k)$ we have \[\ker(\lambda_{m_1})\cap\cdots\cap\ker(\lambda_{m_k})=\ker(\lambda_g).\]
\el
\bdem Since $g=\gcd(m_1,\ldots,m_k)$, for every  $i\leq k$ there is an integer $d_i$ such that $m_i=d_i\cdot g$. Moreover, there are integers $\alpha_1,\ldots,\alpha_k$ such that $g=\alpha_1\cdot m_1+\cdots+\alpha_k\cdot m_k$.\\ 

For one inclusion, notice that if $x\in\ker(\lambda_g)$ then we have $g\cdot x=0$, which implies \[m_i\cdot x=(d_i\cdot g)\cdot x=m_i\cdot 0=0.\] Hence, $x\in \ker(\lambda_{m_1})\cap\cdots \cap \ker(\lambda_{m_k})$. For the other inclusion, notice that if $x\in \ker(\lambda_{m_1})\cap \cdots\cap \ker(\lambda_{m_k})$ then we have
\begin{align*}
g\cdot x&=(\alpha_1\cdot m_1+\cdots+\alpha_k\cdot m_k)\cdot x=\alpha_1(m_1\cdot x)+\cdots+\alpha_k(m_k\cdot x)\\
&=\alpha_1\cdot 0+\cdots+\alpha_k\cdot 0=0.
\end{align*}
Therefore, $x\in \ker(\lambda_g)$.
\edem

\bl \label{Lemma-Coprime} Let $A$ be an abelian group, and let $p,q\geq 1$ with $\gcd(p,q)=1$. Then, \[A/(pqA) \cong (A/pA)\oplus(A/qA).\]
\el
\bdem Consider the canonical map $f:A\to (A/pA)\oplus(A/qA)$ given by $f(x)=(x+pA,x+qA)$. It is clear that $f$ is a group homomorphism, so it is enough to show that $f$ is surjective and that $\ker(f)=pqA$.\\

Since $\gcd(p,q)=1$, there are integers $r,s$ such that $r\cdot p+s\cdot q=1$.

Clearly, if $x\in pqA$, then $x\in pA$ and $x\in qA$, which implies that $f(x)=(pA,qA)=(0,0)$. Conversely, suppose that $x\in \ker(f)$. Then $f(x)=(0,0)$, which implies that $x\in pA$ and $x\in qA$. So, there are elements $y,z\in A$ such that $x=py=qz$. Hence, 

\begin{align*}
r\cdot x=r\cdot (p\cdot y)&=r\cdot (q\cdot z)\\
(1-s\cdot q)\cdot y&=r\cdot q\cdot z\\
y&=s\cdot q\cdot y+r\cdot q\cdot z=q(s\cdot y+r\cdot z).
\end{align*}
By taking $w=s\cdot y+r\cdot z$, we have $y=q\cdot w$, and so $x=p\cdot y=p\cdot (q\cdot w)$. Hence, $x\in pqA$.\\

Let us now see that $f$ is surjective, that is, given $a,b\in A$ there is $x\in A$ such that $f(x)=(a+pA,b+qA)$. Let us take $x=rp\cdot b+sq\cdot a$. We have the following:
\begin{align*}
x&=rp\cdot b+sq\cdot a=rp\cdot b+(1-rp)\cdot a=a+p(r\cdot b-r\cdot a) & &\Rightarrow & &x-a\in pA,\\
x&=rp\cdot b+sq\cdot a=(1-sq)\cdot b+sq\cdot a=b+q(-s\cdot b+s\cdot a) & &\Rightarrow & &x-b\in qA.
\end{align*}
Hence, $f(x)=(a+pA,b+qA)$, as desired.
\edem

\bl \label{Lemma-CountingKernel} Let $A$ be a finite abelian group and $m\geq 1$, with prime decomposition $m=p_1^{t_1}\cdots p_k^{t_k}$. Then,
\[|\operatorname{ker}(\lambda_m)|=\prod_{j=1}^k\left(\prod_{i=0}^{t_j-1} |p_j^{i}A/p_j^{i+1}A|\right).\]
\el
\bdem If $m=p_1^{t_1}\cdots p_k^{t_k}$ is the prime decomposition of $m$, then by Lemma \ref{Lemma-Coprime} 
\begin{align*}
A/mA&\cong A/(p_1^{t_1}\cdots p_k^{t_k}A)\cong (A/p_1^{t_1}A) \oplus \cdots \oplus (A/p_k^{t_k}A).
\end{align*}
Also, notice that for every prime $p$ we have 
\begin{align*}
|A/p^tA|&=\dfrac{|A|}{|p^tA|}=\dfrac{|A|}{|pA|}\cdot\dfrac{|pA|}{|p^2A|}\cdots\dfrac{|p^{t-1}A|}{|p^tA|}=\prod_{i=0}^{t-1}|p^iA/p^{i+1}A|.
\end{align*}
Furthermore, by the fundamental homomorphism theorem, we have $A/\operatorname{ker}(\lambda_m)\cong \operatorname{Im}(\lambda_m)$. Hence, for a finite abelian group $A$ we have $|A|=|\ker(\lambda_m)|\cdot |\operatorname{Im}(\lambda_m)|$, which implies
\begin{align*}
|\ker(\lambda_m)|&=\dfrac{|A|}{|\operatorname{Im}(\lambda_m)|}=|A/mA|=\left|(A/p_1^{t_1}A) \oplus \cdots \oplus (A/p_k^{t_k}A)\right|\\
&=\prod_{j=1}^k|A/p_j^{t_j}A|=\prod_{j=1}^k\left(\prod_{i=0}^{t_j-1}|p_j^{i}A/p_j^{i+1}A|\right). \qedhere
\end{align*}
\edem

\bl \label{Lemma-CountingDivisibility} Let $A$ be a finite abelian group and let $p$ prime, $n\geq 1$ and $r\in \mathbb{Z}$. If $r=p^s\cdot r'$ with $\gcd(r',p)=1$, then the subgroup $H_{p,n,r}=\{x\in A:r\cdot x\in p^nA\}$ has cardinality 
\[|H_{p,n,r}|=\prod_{i=s}^{n-1}|p^iA/p^{i+1}A|\]
\el

\bdem Consider the map $\widehat{\lambda}_r:A\to p^nA$ defined by $\widehat{\lambda}_r(x)=r\cdot x+p^nA$. Notice that this is a homomorphism, and its kernel is precisely the subgroup $H=H_{p,n,r}$. Also, notice that the image of $\widehat{\lambda}_r$ is $rA+p^nA$, which corresponds to is $p^sA$. Hence, using the same idea as the Lemma \ref{Lemma-CountingKernel}
\[|H_{p,n,r}|=|\operatorname{ker}(\widehat{\lambda}_r)|=\dfrac{|p^nA|}{|p^sA|}=\dfrac{|A/p^nA|}{|A/p^sA|}=\dfrac{\prod_{i=0}^{n-1}|p^iA/p^{i+1}A|}{\prod_{i=0}^{s-1}|p^iA/p^{i+1}A|}=\prod_{i=s}^{n-1}|p^iA/p^{i+1}A|. \qedhere\]
\edem

\brmk Notice a feature for \emph{finite} abelian groups that we have used repeatedly in the previous results. If $A$ is an abelian group, and $\lambda:A\to A$ is a definable homomorphism, then by the fundamental theorem of isomorphisms we have $A/\operatorname{ker}(\lambda)\cong \operatorname{Im}(\lambda)$. So, in a finite abelian group, whenever $\ker(\lambda)\neq 0$ we would have that $\operatorname{Im}(\lambda)$ is a \emph{proper} subgroup of $A$. This is false in the general setting. For instance, if we consider the group $A=\mathbb{Z}(p^\infty)^\omega$ (infinitely many copies of the Prüfer $p$-group), then we have $\operatorname{ker}(\lambda_p)\cong (\mathbb{Z}/p\mathbb{Z})^\omega$ is infinite, but $\operatorname{Im}(\lambda_p)=pA=A$ because $A$ is divisible. 
\ermk

\begin{notation}\label{Notation-Groups} Since we will be working throughout the paper with several examples of abelian groups, we would like to set the following conventions:
\bitem
\item $\mathbb{Z}(p^{n})$ will represent the cyclic group of order $p^n$.
\item $\mathbb{Z}(p^{\infty})$ represents Pr\"ufer $p$-group.
\item $\mathbb{Z}_{(p)}$ for the localization of the integers $\mathbb{Z}$ at the prime $p$, i.e., the additive subgroup of the rational numbers with a denominator that is not divisible by $p$.
\item $\mathbb{Q}$ represents the additive group of the rational numbers.
\eitem
\end{notation}
\subsection{Model theory of abelian groups}
\label{Subsec:ModelTheoryAG}
Our eventual goal is to characterize exactly which abelian groups arise as ultraproducts of finite abelian groups. The route to this classification runs through Szmielew’s elimination theory: every first-order property of an abelian group is controlled by a countable family of numerical invariants attached to its positive primitive definable subgroups. We build this machinery in two steps. First, we introduce Pr\"ufer’s Lemma which reduces the complexity of positive primitive formulas (Lemma \ref{PruferPPF}) as an application of the Smith Normal Form (Fact \ref{Fact:SmithNormalForm}).\\

Then Szmielew’s theorem (Theorem \ref{SzmielewTheorem}) lets us understand an arbitrary abelian group by a canonical representative, called a Szmielew group, on which every one of these invariants can be computed explicitly, which is later used in the classification results of Section \ref{sect:Characterizationpsfgroups}. We start with the following fact.

\bh [Smith Normal Form (see \cite{Smith1861})] \label{Fact:SmithNormalForm} Let $T$ be a matrix over $\mathbb{Z}$ of rank $\nu$. There are invertible matrices $U,V$ with integer coefficients such that $UTV=D$ is a matrix of the form

\[
D=\begin{bmatrix}
  \begin{matrix}
  c_1 & 0 & \cdots & 0 \\
  0 & c_2 & \cdots & 0 \\
  \vdots & \vdots & \ddots & \vdots\\
  0 & 0 &\cdots & c_{\nu}
  \end{matrix}
  & \rvline & & \bigzero & \\
\hline
& \rvline &\\
  \bigzero & \rvline & & \bigzero
\end{bmatrix}
\]with integers $c_1,c_2,\cdots,c_\nu\geq 1$
\eh
In abelian groups, we will show that it is possible to simplify formulas to resolutions of systems of equations. We call those basic systems \emph{positive primitive formulas} (or simply \emph{p.p. formulas}), with the form: \[\exists \ov{y}(\psi_1(\ov{x};\ov{y})\wedge \psi_2(\ov{x};\ov{y})\wedge \cdots \wedge \psi_k(\ov{x};\ov{y}))\]
where each $\psi_i(\ov{x};\ov{y})$ is a linear equation on its corresponding variables.\\

Furthermore, the positive primitive formulas define a subgroup of any abelian group as they are solutions to a matrix equation. In particular, using the Smith Normal Form we can encode these matrices in the following way:
\bl [Prüfer] \label{PruferPPF} Every positive primitive formula $\varphi(\ov{x})$ is equivalent, modulo the theory of abelian groups, to a finite conjunction of formulas of the form $t(\ov{x})=0$ or $p^n|t(\ov{x})$ where $t(\ov{x})$ is a term, $p$ is a prime and $n$ is a positive natural number.
\el
\bdem It is easy to check that in the language $\Le=\{+,0,-\}$ of abelian groups every term $t(x_1,\ldots,x_n)$ is equivalent modulo the theory of abelian groups to an expression of the form $k_1\cdot x_1+\cdots+k_n\cdot x_n$ for some integers $k_1,\ldots,k_n\in\mathbb{Z}$. Let $\varphi(\ov{x})$ be a positive primitive formula with $|\ov{x}|=n$, and let us write it as $\exists \ov{y}(R\ov{x}+Q\ov{y})$ where $R$ and $Q$ are matrices of sizes $k\times n$ and $k\times m$ (respectively) with integer coefficients. Using the Smith Normal Form, there are invertible matrices $U,V$ with integer coefficients such that $UQV$ is diagonal. Hence, we have
\begin{align*}
R\ov{x}+Q\ov{y}=0&\Leftrightarrow UR\ov{x}+UQ\ov{y}=0\Leftrightarrow UR\ov{x}+UQVV^{-1}\ov{y}=0\Leftrightarrow UR\ov{x}=UQV(-V^{-1}\ov{y})
\end{align*}
By rewriting $\ov{z}=-V^{-1}\ov{y}$, the formula $\varphi(\ov{x})$ is equivalent to the existence of a solution $\ov{z}$ of the matrix equation $UR\ov{x}=UQV\ov{z}$. Furthermore, since $UQV$ is a diagonal matrix and $UR\ov{x}$ defines a system of linear combinations of the variables in $\ov{x}$, we can rewrite $\varphi(\ov{x})$ as
\[\exists \ov{z}(t_1(\ov{x})=\lambda_1z_1\wedge \cdots \wedge t_k(\ov{x})=\lambda_kz_k),\]
where $\lambda_i$ are the diagonal elements of the matrix $UQV$. Moreover, notice that each of the conjunctions contains at most one of the variables $z_i$, so this is equivalent to
\[\exists z_1(t_1(\ov{x})=\lambda_1z_1) \wedge \exists z_2(t_2(\ov{x})=\lambda_2z_2)\wedge \cdots \wedge \exists z_m(t_m(\ov{x};\ov{y})=\lambda_mz_m).\]

We now study each of these conjunctions. If $\lambda_i=0$, we have a formula of the form $t_i(\ov{x})=0$. Otherwise, we have the formula $\exists z_i(t_i(\ov{x})=\lambda_iz_i)$, which is equivalent to state that $\lambda_i|t_i(\ov{x})$. By decomposing $\lambda_i$ into prime factors, this is equivalent to a finite conjunction of formulas of the form $p^n|t_i(x)$ for some primes $p$ and $n>0$. Hence, the primitive positive formula $\varphi(\ov{x})$ is equivalent to a finite conjunction of formulas of the form $t(\ov{x})=0$ or $p^n|t(\ov{x})$, as desired.
\edem

Two direct consequences of Lemma \ref{PruferPPF} are the following:
\bl \label{ConmutacionPPFormulaProducto}Let $\langle A_i:i\in I\rangle$ be a non-empty family of abelian groups and $\varphi(\ov{x})$ a primitive positive formula. Then
\[\varphi\left(\prod_{i\in I}A_i\right)=\prod_{i\in I}\varphi(A_i).\]
\el

\bl \label{ConmutacionCocienteProducto} Let $\langle A_i:i\in I\rangle$ be a non-empty family of abelian groups, and let $\langle M_i:i\in I\rangle,\langle N_i:i\in I\rangle$ be two non-empty families of subgroups satisfying $N_i\leq M_i\leq A_i$ for every $i\in I$. Then,
\[\prod_{i\in I}M_i\Big/\prod_{i\in I}N_i\cong \prod_{i\in I}M_i/N_i\ \ .\]
\el

Now we want to replace an abelian group by an explicit family of abelian groups which have an easier invariant computation. This is a family that we will obtain via Szmielew's theorem (Theorem \ref{SzmielewTheorem}).\\

The collection of Szmielew groups are obtained from four families of building blocks: Finite cyclic groups, Pr\"ufer groups, localizations of $\Z$ and the group of rationals $\Q$. We will summarize the invariants of each building block in Table \ref{TablaInvariantesSzmielew} since these computations are the base of Section \ref{sect:Characterizationpsfgroups}.
\bd[Szmielew group] \label{SzmielewGroupDefinition} A \textbf{Szmielew group} is a finite or countable abelian group of the form:
\begin{align*}
    \bigoplus_{p \text{ prime}} \left[ \bigoplus_{n>0} \mathbb{Z}(p^{n})^{\kappa_{(p,n-1)}} \oplus \mathbb{Z}(p^{\infty})^{\lambda_{p}} \oplus \mathbb{Z}^{\mu_p}_{(p)} \right] \oplus \mathbb{Q}^{\nu},
\end{align*}
where $\kappa_{(p,n-1)},\lambda_{p},\mu_{p},\nu$ are cardinals less than or equal to $\omega$.
\ed

\bt[Szmielew's theorem] \label{SzmielewTheorem} Every abelian group is elementary equivalent to a Szmielew group. (See \cite{Szmielew1955})
\et
\bd[Szmielew invariants] (See \cite{Szmielew1955})\label{SzmielewInvariantDefinition} For any abelian group $A$, $n\in\mathbb{N}$ and prime $p$, we write $A[p]$ for the subgroup $\{a\in A : pa=0\}$, and denote $p^nA[p]$ the realizations of the positive primitive formula $\varphi(x):=px=0\hspace{0.1cm}\land\hspace{0.1cm}p^n|x$. We define the \textbf{Szmielew invariants} of $A$ as the following numerical invariants:
\begin{align*}
    &\operatorname{U}(p,n;A) = \left|p^n A[p]/p^{n+1}A[p]\right|, & & &\operatorname{D}(p,n;A) = \left|p^n A[p]\right|,\\
    &\operatorname{Tf}(p,n;A) = \left|p^n A/p^{n+1}A\right|, & & &\operatorname{Exp}(p,n;A) = \left|p^n A\right|.
\end{align*}
Here, $U$ is a refers to the Ulm invariant (as in \cite{fuchs1970infinitegroups} p.154), $D$ is the divisible invariant, $Tf$ represents the torsion free invariant and $Exp$ is the exponent invariant.
\ed
Note that all the Szmielew invariants are representing the size of a quotient of definable sets by a positive primitive formula.
\bo Given an invariant $F$ where the positive primitive formulas $\varphi$ and $\psi$ define the numerator and denominator respectively, it is possible to write a sentence expressing that the invariant $F=\text{Inv}\left(\varphi,\psi,A\right)$ has a value at least $k$ (for $k<\omega$) via the following sentence $\hat{\sigma}_{F,k}$:
\begin{align*}
    \hat{\sigma}_{F,k}:=\exists z_{1},\dots,z_{k} \left( \bigwedge\limits_{i=1}^{k}(\phi(z_{i})) \land \bigwedge\limits_{1\leq i < j\leq k} \neg (\psi(z_{i}-z_{j}))\right)
\end{align*}
Moreover, the following sentence expresses that the invariant $F$ has value $k$:
\begin{align*}
    \sigma_{(F,k)} : \hat{\sigma}_{(F,k)} \land \neg \hat{\sigma}_{(F,k+1)}.
\end{align*}
These sentences will be called \textit{Szmielew invariant sentences} denoted by:
\begin{align*}
    \sigma_{p,n,k}^{U}, \sigma_{p,n,k}^{D}, \sigma_{p,n,k}^{Tf}, \sigma_{p,n,k}^{Exp}
\end{align*}
for a prime $p$ and $n\in\mathbb{N}$.
\eo
\bt\label{thm:SzmiliewInvariantsEquiv} For every abelian group $A$, the theory $\operatorname{Th}(A)$ is determined by the Szmielew invariants of $A$. Equivalently, every first-order sentence in the language of abelian groups is equivalent modulo $T_{ab}$ to a boolean combination of Szmielew invariant sentences.
\et
\bc Let $A$, $B$ be two abelian groups. $A\equiv B$ if and only if $A$ and $B$ have the same Szmielew invariants.
\ec
Observe that Theorems \ref{SzmielewTheorem} and \ref{thm:SzmiliewInvariantsEquiv} describe that studying Szmielew invariants of abelian groups requires an understanding of how these numerical invariants behave in direct sums and products of the structures. The following lemma addresses that inquiry.
\bl Let $\langle A_i:i\in I\rangle$ be a non-empty family of abelian groups. Then, for every Szmielew invariant $F$ we have
\[F\left(\prod_{i\in I}A_i\right)=\sup\left\{\prod_{i\in J}F(A_i):J\subseteq^{\text{fin}}I\right\}=F\left(\bigoplus_{i\in I}A_i\right).\]
\label{BaurInvariantsForSums}
\el
\bdem Let $\varphi(\ov{x}),\psi(\ov{x})$ be two positive primitive formulas for $\Le$. By Lemma \ref{ConmutacionPPFormulaProducto} we have $\varphi\left(\prod_{i\in I}A_i\right)=\prod_{i\in I}\varphi(A_i)$, and the same holds for $\psi$. In particular, we can write $F$ in terms of two positive primitive formulas $\varphi$ and $\psi$ such that $F(A)=F_{\varphi,\psi}(A)=|\varphi(A)/(\varphi(A)\cap\psi(A))|$ for every abelian group $A$, we have by Lemma \ref{ConmutacionCocienteProducto} that
\begin{align*}
F\left(\prod_{i\in I}A_i\right)&=\left|\varphi\left(\prod_{i\in I}A_i\right)\huge/(\varphi\wedge\psi)\left(\prod_{i\in I}A_i\right)\right|=\left|\prod_{i\in I}\left(\varphi(A_i)\Big/(\varphi\wedge\psi)(A_i)\right)\right|.
\end{align*}
Note also that the Szmielew invariants are always $\geq 1$, so we also have for every $J\subseteq^{\text{fin}}I$ that $\prod_{i\in J}F(A_i)\leq \prod_{i\in I}F(A_i)$, so
\[F\left(\prod_{i\in I}A_i\right)\geq\sup\left\{\prod_{i\in J}F(A_i):J\subseteq^{\text{fin}}I\right\}.\]
We now show that the equality holds, arguing by cases on the value of $F\left(\prod_{i\in I}A_i\right)$. If it is equal to $1$, the equality follows immediately as that is the smallest value that a Szmielew invariant can take.\\

Suppose now that $F\left(\prod_{i\in I}A_i\right)=n<\infty$. Then there are $n$ distinct representatives $\alpha_1,\ldots,\alpha_n$ satisfying that $\alpha_i-\alpha_j\not\in (\varphi\wedge\psi)\left(\prod_{i\in I}A_i\right)$ whenever $1\leq i<j\leq n$. So, for each pair $i<j$ there is an index $k_{ij}\in I$ such that the $k_{ij}$-th coordinates of $\alpha_i$ and $\alpha_j$ are different modulo $(\varphi\wedge\psi)(A_{k_{ij}})$. Note that $J_0=\{k_{ij}:1\leq i<j\leq n\}$ is a finite subset of $I$, and there are $n$ different representatives $\alpha_1,\ldots,\alpha_n$ of $\varphi\left(\prod_{i\in J_0}A_i\right)$ modulo $(\varphi\wedge\psi)\left(\prod_{i\in J_0}A_i\right)$, so $F\left(\prod_{i\in J_0}A_i\right)=n=F\left(\prod_{i\in I}A_i\right)$. Hence, by the bound we proved before, the equality holds.\\

Finally, suppose that $F\left(\prod_{i\in I}A_i\right)=\infty$. This implies that for every $n$ it is possible to find $n$ distinct representatives in $\varphi\left(\prod_{i\in I}A_i\right)$ modulo $(\varphi\wedge\psi)\left(\prod_{i\in I}A_i\right)$. By the same argument as above, there are always finite subsets $I_n\subseteq I$ such that $\varphi(\prod_{i\in I_n}A_i)$ contains at least $n$ different representatives modulo $(\varphi\wedge\psi)\left(\prod_{i\in I_n}A_i\right)$, and so $\prod_{i\in I_n}F(A_i)\geq n$. Therefore, 
\[\sup\left\{\prod_{i\in J}F(A_i):J\subseteq^{\text{fin}}I\right\}=\sup\left\{\prod_{i\in I_n}F(A_i):n<\omega\right\}=\infty=F\left(\prod_{i\in I}A_i\right).\]

Note also that this equality implies that $F\left(\prod_{i\in I}A_i\right)=F\left(\bigoplus_{i\in I}A_i\right)$, and this finishes the proof.
\edem
\bc \label{InvariantesSumaGrupos} Let $\langle A_i:i\in I\rangle$ be a finite non-empty family of abelian groups. Then, for every Szmielew invariant $F$, 
\[F\left(\prod_{i\in I}A_i\right)=\prod_{i\in I}F(A_i)=F\left(\bigoplus_{i\in I}A_i\right).\]
\ec

Furthermore, Table \ref{TablaInvariantesSzmielew} summarizes the values of the Szmielew invariants in each of the basic factors of Szmielew groups for $p,q$ prime numbers, $k>0$ and $n\geq0$.\\

\begin{table}[h!]
\centering
\begin{tabular}{|c|c|c|c|c|}
\hline
\textbf{Invariant} & $\mathbb{Z}(p^{k})$                                                                                & $\mathbb{Z}(p^{\infty})$                                                                & $\mathbb{Z}_{(p)}$                                                              & $\mathbb{Q}$        \\ \hline
$\operatorname{U}(q,n;-)$          & \begin{tabular}[c]{@{}c@{}}$p$ if $(q=p,n=k-1)$\\ $1$ (otherwise)\end{tabular}         & 1                                                                         & 1                                                                         & 1        \\ \hline
$\operatorname{D}(q,n;-)$          & \begin{tabular}[c]{@{}c@{}}$p$ if $(q=p,n<k)$\\ $1$ (otherwise)\end{tabular}           & \begin{tabular}[c]{@{}c@{}}$p$ if $(q=p)$\\ $1$ (otherwise)\end{tabular} & 1                                                                         & 1        \\ \hline
$\operatorname{Tf}(q,n;-)$         & \begin{tabular}[c]{@{}c@{}}$p$ if $(q=p,n<k)$\\ $1$ (otherwise)\end{tabular}           & 1                                                                         & \begin{tabular}[c]{@{}c@{}}$p$ if $(q=p)$\\ $1$ (otherwise)\end{tabular} & 1        \\ \hline
$\operatorname{Exp}(q,n;-)$        & \begin{tabular}[c]{@{}c@{}}$p^{k-n}$ if $(q=p,n<k)$\\ $p^{k}$ (otherwise)\end{tabular} & $\infty$                                                                  & $\infty$                                                                  & $\infty$ \\ \hline
\end{tabular}
\caption{Szmielew invariant of the basic factors of Szmielew groups.}
\label{TablaInvariantesSzmielew}
\end{table}

In this article, we will compute the invariants of the finite cyclic completely, and characterize the Torsion free invariant for all Szmielew groups.

\bej\label{ExampleCyclicGroups} Let $p$ be a prime and $k>0$, consider the finite abelian cyclic group of order $p^k$ which we denote as $(\Z(p^k),+)$. First observe that for any prime $q\neq p$, $\Z(p^k)$ does not have $q$-torsion and so $\Z(p^k)[q]=\{0\}$, and by $q$-divisibility in the group, we get that $q^n\Z(p^k) = \Z(p^k)$ for any $n$. Obtaining the following invariants:
\begin{align*}
    \operatorname{U}(q,n;\mathbb{Z}(p^{k})) &= \left|q^n \mathbb{Z}(p^{k})[q]/q^{n+1}\mathbb{Z}(p^{k})[q]\right| = \left|\{0\}/\{0\}\right|=1\\
    \operatorname{D}(q,n;\mathbb{Z}(p^{k})) &= \left|q^n \mathbb{Z}(p^{k})[q]\right|=\left|\{0\}\right|=1\\
    \operatorname{Tf}(q,n;\mathbb{Z}(p^{k})) &= \left|q^n \mathbb{Z}(p^{k})/q^{n+1}\mathbb{Z}(p^{k})\right|=\left|\mathbb{Z}(p^{k})/\mathbb{Z}(p^{k})\right|=1\\
    \operatorname{Exp}(q,n;\mathbb{Z}(p^{k})) &= \left|q^n \mathbb{Z}(p^{k})\right| =|\mathbb{Z}(p^{k})|=p^{k}
\end{align*}
Now, observe that an element $x\in\Z(p^k)[p]$ must be of the form $x=p^{k-1}\cdot r$ for some $r$ between $0$ and $p-1$ which can be expressed in the following form:
\begin{align*}
    \Z(p^{k})[p]=\{p^{k-1}\cdot r : r\in\{0,1,\cdots,p-1\}\}.
\end{align*}
Thus for the case when $n\geq k$, $p^n\Z(p^k)[p]={0}$, otherwise, $p^n\Z(p^k)[p]=\Z(p^k)[p]$. On the other hand, for any $y\in p^{n}\Z(p^{k})$, it has the form $y=z\cdot p^n$ for some $0\leq z\leq p^{k-n}$ since $0\leq y<p^k$. Then the subgroup $p^{n}\mathbb{Z}(p^{k})=\{z\cdot p^{n} : 0\leq z<p^{k-n}\}$ and so, if $n\geq k$, then $p^{n}\Z(p^{k})=\{0\}$. Otherwise $|p^{n}\Z(p^{k})|=p^{k-n}$.\\

Obtaining the following invariants:
\begin{align*}
    \operatorname{U}(p,n;\mathbb{Z}(p^{k})) &= \left|p^n \mathbb{Z}(p^{k})[p]/p^{n+1}\mathbb{Z}(p^{k})[p]\right|= \frac{|p^n \mathbb{Z}(p^{k})[p]|}{|p^{n+1} \mathbb{Z}(p^{k})[p]|}=\left\{ \begin{array}{lc}
             p &   \text{if }  n=k-1 \\
             \\ 1 &  \text{otherwise} 
             \end{array}
   \right.\\
    \operatorname{D}(p,n;\mathbb{Z}(p^{k})) &= \left|p^n \mathbb{Z}(p^{k})[p]\right|=\left\{ \begin{array}{lc}
             p &   \text{if }  n<k \\
             \\ 1 &  \text{otherwise} 
             \end{array}
   \right.\\
    \operatorname{Tf}(p,n;\mathbb{Z}(p^{k})) &= \left|p^n \mathbb{Z}(p^{k})/^{n+1}\mathbb{Z}(p^{k})\right|= \frac{|p^n \mathbb{Z}(p^{k})|}{|p^{n+1} \mathbb{Z}(p^{k})|}=\left\{ \begin{array}{lc}
             p &   \text{if }  n<k \\
             \\ 1 &  \text{otherwise} 
             \end{array}
   \right.\\
    \operatorname{Exp}(p,n;\mathbb{Z}(p^{k})) &= \left|p^n \mathbb{Z}(p^{k})\right| =\left\{ \begin{array}{lc}
             p^{k-n} &   \text{if }  n<k \\
             \\ 1 &  \text{otherwise} 
             \end{array}
   \right.\\
\end{align*}
\eej
\bej Now, we will compute the Torsion free invariant $\operatorname{Tf}(p,n;\cdot)$ for $p$ prime and $n\geq0$ for $\Q$, $\Z(p^\infty)$ and $\Z_{(p)}$ in that same order. Note that both $(\Q,+)$, $\Z(p^\infty)$ are divisible abelian groups, so for any combination of $p$ and $n$, the $p^n$-divisible component of each group is everything. Hence, $\operatorname{Tf}(p,n;A) = \left|p^n A/p^{n+1}A\right|=|A/A|=1$, for $A=\Q$ or $A=\Z(p^\infty)$.\\

Finally, the case of $\Z_{(p)}$ is the only non-trivial computation for this example, since it is not divisible. However, for any prime $q\neq p$, we have that $q^n\Z_{(p)}=\Z_{(p)}$ and $\operatorname{Tf}(q,n;\Z_{(p)})=1$. Then we characterize the subgroup of divisible elements by $p^n$ explicitly by:
\begin{align*}
    p^n\mathbb{Z}_{(p)}=\left\{x\in\mathbb{Z}_{(p)}: x=\dfrac{p^na}{b},\hspace{0.2cm}(p^na, b)=1,\hspace{0.2cm}p\nmid b\right\}.
\end{align*}
Moreover, any element $y\in p^n\Z_{(p)}$ has the form $\dfrac{p^n\cdot a}{b}$ where $(p^n\cdot a,b)=1$ and $p\nmid b$. As $a$ is an integer, it can be written in the form $t\cdot p+r$ where $r$ is an integer between $0$ and $p-1$, and so:
\begin{align*}
    y&=\dfrac{p^n \cdot a}{b} = p^n \dfrac{t\cdot p + r}{b} = p^n\left(\dfrac{t\cdot p}{b} + \dfrac{r}{b}\right)=\dfrac{t\cdot p^{n+1}}{b} + \dfrac{r\cdot p^n}{b}.
\end{align*}
Thus any element of $p^n\Z_{(p)}$ is equivalent to an element of the form $\dfrac{r\cdot p^n}{b}$ modulo $p^{n+1}\Z_{(p)}$. So, we conclude that the quotient between both subgroups is of the form:
\begin{align*}
    p^n\mathbb{Z}_{(p)}/p^{n+1}\mathbb{Z}_{(p)}=\left\{\left[0\right],\left[\dfrac{p^n}{b}\right],\left[\dfrac{2p^n}{b}\right]\dots\left[\dfrac{(p-1)p^n}{b}\right]\right\}.
\end{align*}
which has exactly $p$ different equivalence classes for any $n$, equivalently, $\operatorname{Tf}(p,n;\Z_{(p)})=p$.
\eej
\section{Characterization of pseudofinite abelian groups}\label{sect:Characterizationpsfgroups}
This section establishes our main classification results for pseudofinite abelian groups, in three steps. In Subsection \ref{Subsec:PsfDirectSums} we show that the class of pseudofinite abelian groups is closed under countable direct sums, which reduces the analysis of Szmielew groups to the combination of exponents that generate pseudofinite factors.\\

Later, in Subsection \ref{Subsec:PsfAG} we prove our main theorem (Theorem \ref{Theorem:ClassificationPsfAbelianGroups}), a complete characterization of pseudofinite abelian groups in terms of their Szmielew invariants. Finally, in Subsection \ref{Subsec:PsfWStableAG} we adapt this characterization to the $\omega$-stable case (Theorem \ref{ClassificationWStablePsfAbGroups}) by combining Theorem \ref{Theorem:ClassificationPsfAbelianGroups} with Macintyre's theorem (Theorem \ref{TeoremaMacintyre2.0}).

\subsection{Preservation of pseudofiniteness under direct sums}
\label{Subsec:PsfDirectSums}
We begin with the ultrafilter machinery needed to show pseudofiniteness is preserved under direct sums. The key technical obstacle is that a direct sum of two ultraproducts is not an ultraproduct in a direct way. However, Lemmas \ref{FilterForDirectSum} and \ref{SumTwoUltraproducts} construct, from two ultrafilters $\mathcal{U}$ and $\mathcal{V}$ on two index sets, a single ultrafilter $\mathcal{W}$ on their product witnessing the expected product construction of two ultraproducts, which immediately gives us closure of pseudofinite structures under finite direct sums (Theorem \ref{SumaDirectaDosGAPseudo}).\\

Extending this to countable direct sums takes two further ingredients: An argument identifying a countable direct sum with an ultraproduct of finite partial sums (Theorem \ref{EquivUltraSumaEnum}), and a general fact that an ultraproduct of ultraproducts of finite structures is again an ultraproduct of finite structures. (Lemma \ref{ProductUltrafilter}, and Theorem \ref{Ultraproductodeultraproductosestfinitas}) Together we obtain Corollary \ref{CountableSumPseudofinite} which concludes that pseudofiniteness is preserved under countable direct sums. An independent proof of the closure of pseudofinite structures under products was constructed by D'Aquino and Macintyre in \cite{DAquino-Macintyre2026}.\\

We start the analysis with the following two lemmas, which have fairly straightforward proofs.

\bl \label{FilterForDirectSum} Let $\U$ and $\mathcal{V}$ be ultrafilters over the index sets $I$ and $J$ respectively. There is an ultrafilter $\mathcal{W}$ over $I\times J$ extending $\mathcal{U}\otimes \mathcal{V}$ where the latter is defined as:
\begin{align*}
    \mathcal{U}\otimes \mathcal{V} = \{ U\times V : U\in\mathcal{U} \text{ and } V\in\mathcal{V}\}.
\end{align*}
\el

\bl \label{SumTwoUltraproducts} Let $\mathcal{U},\mathcal{V}$ be two ultrafilters on index sets $I,J$. Let $\langle \mathcal{M}_i\rangle$, $\langle\mathcal{N}_{j} : j\in J\rangle$ be two non-empty families of abelian groups. Then, there is an ultrafilter  $\mathcal{W}$ over $I\times J$ such that
\begin{align*}
    \prod_{\mathcal{U}}\mathcal{M}_{i} \oplus \prod_{\mathcal{V}}\mathcal{N}_{j} \equiv \prod_{\mathcal{W}}  (\mathcal{M}_{i} \oplus \mathcal{N}_{j})\hspace{0.5cm}.
\end{align*}
\el

\bt \label{SumaDirectaDosGAPseudo} If $G$ and $H$ are two pseudofinite abelian groups, then $G\oplus H$ is pseudofinite.
\et
\bdem Since $G$ and $H$ are pseudofinite, there are two ultrafilters $\mathcal{U}$ and $\mathcal{V}$ over sets of indices $I,J$, and two families of finite abelian groups $\langle\mathcal{G}_{i} : i\in I\rangle$, $\langle\mathcal{H}_{j} : j\in J\rangle$ such that:
\begin{align*}
    G \equiv \prod_{\mathcal{U}} \mathcal{G}_{i} \text{\ \ and\ \  }    H \equiv \prod_{\mathcal{V}} \mathcal{H}_{i}.
\end{align*}
By Lemma \ref{SumTwoUltraproducts}, there is an ultrafilter $\mathcal{W}$ extending $\mathcal{U}\otimes\mathcal{V}$, which satisfies
\[\prod_{\mathcal{U}}\mathcal{G}_{i} \oplus \prod_{\mathcal{V}}\mathcal{H}_{j} \equiv \prod_{\mathcal{W}}  (\mathcal{G}_{i} \oplus \mathcal{H}_{j}).\]

Since $G \equiv \prod\limits_{\mathcal{U}} \mathcal{G}_{i}$ and $H \equiv \prod\limits_{\mathcal{V}} \mathcal{H}_{j}$, for every Szmielew invariant $F$ we have  $F(G) = F\left(\prod\limits_{\mathcal{U}} \mathcal{G}_{i}\right)$ and $F(H) = F\left(\prod\limits_{\mathcal{V}} \mathcal{H}_{j}\right)$. Also, by Lemma \ref{SumTwoUltraproducts}, we have
\[F\left(\prod\limits_{\mathcal{U}}\mathcal{G}_{i} \oplus \prod\limits_{\mathcal{V}}\mathcal{H}_{j}\right)= F\left(\prod\limits_{\mathcal{W}}  (\mathcal{G}_{i} \oplus \mathcal{H}_{j})\right).\]
With these equalities, we can calculate the Szmielew invariant $F(G\oplus H)$ using Corollary \ref{InvariantesSumaGrupos}:
\begin{align*}
    F(G\oplus H)&=F(G)\cdot F(H) = F\left(\prod\limits_{\mathcal{U}} \mathcal{G}_{i}\right)\cdot F\left(\prod\limits_{\mathcal{V}} \mathcal{H}_{j}\right)\\ &=F\left(\prod_{\mathcal{U}}\mathcal{G}_{i} \oplus \prod_{\mathcal{V}}\mathcal{H}_{j}\right)=F\left(\prod_{\mathcal{W}}  (\mathcal{G}_{i} \oplus \mathcal{H}_{j})\right).
\end{align*}
Therefore, since the equality holds for every Szmielew invariant, we conclude that $G\oplus H \equiv \prod\limits_{\mathcal{W}}  (\mathcal{G}_{i} \oplus \mathcal{H}_{j})$, and so $G\oplus H$ is pseudofinite.
\edem
An immediate consequence of the previous result is that the direct sum of finitely many pseudofinite abelian groups is pseudofinite. We now show that this actually holds for direct sums with countably many summands.

\bt \label{EquivUltraSumaEnum} Let $\mathcal{U}$ be a non-principal ultrafilter on $\mathbb{N}$, $\langle G_{n}:n\in\mathbb{N}\rangle$ a family of countably many abelian groups and let $G =\bigoplus\limits_{n<\omega} G_{n}$. Define $\mathcal{H}_{n} = G_{1}\oplus \dots \oplus G_{n}$, and consider the group $H=\prod\limits_{\mathcal{U}}\mathcal{H}_{n}$. Then, $G\equiv H$.
\et
\bdem Let $F$ be an arbitrary Szmielew invariant. First, notice that

\[F(G) = F\left( \bigoplus\limits_{n<\omega} G_{n}\right) = \sup\left\{\prod\limits_{j\in J} F(G_{j}) : J\subseteq^{\text{fin}} \mathbb{N}\right\}.\]

Also, for every $n\in\mathbb{N}$ we have
\[F(\mathcal{H}_{n}) = F (G_{1}\oplus \dots \oplus G_{n}) = \sup\left\{\prod\limits_{j\in J} F(G_{j}) : J\subseteq^{\text{fin}} \{1,\dots,n\}\right\}.\]

To show that $G$ and $H$ are elementary equivalent, we will see that $F(G) = F(H)$. We can consider different cases:\\

\underline{\textbf{Case 1: $F(G)=1$:}} In this case we have $F(G_{k})=1$ for every $k<\omega$, and so \[F(\mathcal{H}_n)=F(G_1)\cdot F(G_2)\cdots F(G_n)=1.\] This implies that $\mathcal{H}_{n}$ models the invariant sentence $\hat{\sigma}_{(F,1)}$ for each $n<\omega$, and by \Los' Theorem, $H\models\hat{\sigma}_{(F,1)}$. This implies that $F(H)=1$.\\

\underline{\textbf{Case 2: $F(G)=n<\infty$:}} In this case we have $\displaystyle{n=\sup\left\{\prod\limits_{j\in J} F(G_{j}) : J\subseteq^{fin} \mathbb{N}\right\}.}$

So, there is a finite set $\hat{J}\subseteq^{\text{fin}}\mathbb{N}$ such that $n$ is equal to the product of the invariants $F(G_j)$ with $j\in \widehat{J}$. Since $\hat{J}$ is finite, we can take $m=\max(\hat{J})$.\\

So, for every natural number $k>m$ we have 
\[F(\mathcal{H}_{k}) = \sup\left\{\prod\limits_{j\in J} F(G_{j}) : J\subseteq^{\text{fin}} \{1,\dots,k\}\right\} = \prod\limits_{j\in\hat{J}} F(G_{j}) = n.\]
This implies that $\mathcal{H}_{k}\models\hat{\sigma}_{(F,n)}$ for every $k\geq m$, and so by \Los' Theorem we have $H\models\hat{\sigma}_{(F,n)}$, which implies $F(H)=n$.\\

\underline{\textbf{Case 3: $F(G)=\infty$:}} This implies that $F(G) = \sup\left\{\prod\limits_{j\in J} F(G_{j}) : J\subseteq^{\text{fin}} \mathbb{N}\right\}=\infty$, and so for every $n<\omega$ there is a finite subset $J_n\subseteq^{\text{fin}}\mathbb{N}$ such that $\prod\limits_{j\in J_{n}} F(G_{j})>n$. Since $J_{n}$ is finite, we can take $q_{n}=\max(J_{n})$, and by an argument similar to the one used in Case 2, for every $k\geq q_{n}$ we have $\displaystyle{F(\mathcal{H}_{k}) \geq \prod\limits_{j\in J_{n}} F(G_{j}) >n}$. Hence, $\left\{k\in\mathbb{N}: \mathcal{H}_{k}\models\neg\hat{\sigma}_{(F,n)}\right\} \in \mathcal{U}$, and by \Los' Theorem we have $H\models\neg\hat{\sigma}_{(F,n)}$. Since this holds for every $n<\omega$, we conclude that $F(H)=\infty$.\\

With all these, we conclude that $F(G)=F(H)$ for every Szmielew invariant $F$. Therefore, by Theorem \ref{thm:SzmiliewInvariantsEquiv}, the groups $G$ and $H$ are elementarily equivalent.
\edem
Note that in order to apply Theorem \ref{EquivUltraSumaEnum} to show that countable direct sums of pseudofinite abelian groups are pseudofinite, we need to characterize that an ultraproduct of ultraproducts of finite abelian groups can be reduced to a single ultraproduct of finite abelian groups. This is a fact that holds in general for ultraproducts of $\Le$-structures.

\bl \label{ProductUltrafilter} Let $\mathcal{U}$ be an ultrafilter over an index set $I$, and for each $i\in I$ let $\mathcal{F}_{i}$ be an ultrafilter over an index set $J_i$. Consider now the set of indices
\begin{align*}
    K=\{(i,j):i\in I, j\in J_{i}\},
\end{align*}
and define on $K$ the collection $\mathcal{W}$ of subsets of $K$ satisfying the following property:
\begin{center}
    $X\in\mathcal{W}$ if and only if $\{i\in I: \{j\in J_{i}:(i,j)\in X\}\in\mathcal{F}_{i}\}\in\mathcal{U}$.
\end{center}
Then, $\mathcal{W}$ is an ultrafilter on $K$. 
\el

\bt \label{Ultraproductodeultraproductosestfinitas} Every ultraproduct of ultraproducts of finite structures is isomorphic to an ultraproduct of finite structures. 
\et

\bc\label{CountableSumPseudofinite} Every countable sum of pseudofinite abelian groups is pseudofinite.
\ec
\brmk Notice that the results in Theorems \ref{SumaDirectaDosGAPseudo}, \ref{EquivUltraSumaEnum} and \ref{Ultraproductodeultraproductosestfinitas} only used a careful analysis of Szmielew invariants and their behavior under fundamental operations of groups. By a general version of Szmielew's Theorem to modules (see the Baur-Monk Theorem in \cite{Baur1976}), all these results can be easily adapted to products of finite and pseudofinite $R$-modules for a fixed ring $R$. 
\ermk
\subsection{Pseudofinite abelian groups via Szmielew invariants}
\label{Subsec:PsfAG}
We now start the classification of pseudofinite abelian groups using the Szmielew invariants summarized in Table \ref{TablaInvariantesSzmielew}. We also provide a classification of the Szmielew groups that are elementarily equivalent to infinite ultraproducts of finite abelian groups. This classification appears to be similar to Theorem A.2.14 in \cite{Hodges-ModelTheory}. However, given the difficulty of finding an exact source, we reproduced the results of the theorem with an independent method.

\bej We will prove that both $(\Z_{(p)},+)$ and $(\Z_{p^{\infty}},+)$ are not pseudofinite for $p$ prime by studying the behavior of the function 
\begin{align*}
    f(x):=p\cdot x = \underbrace{x+x+\ldots+x}_{p\text{-times}},
\end{align*}
which is definable by the formula $\theta(a,b):=p\cdot a=b$. Observe that $f$ is an injective function in $(\Z_{(p)},+)$ but it is not surjective as there is no element $x_0\in\Z_{(p)}$ such that $f(x_0)=1$. Since $f(x)$ is definable, then $(\Z_{(p)},+)$ is not pseudofinite.\\

On the other hand, since $\Z_{p^{\infty}} = \Z[1/p]/\Z$, the function $f(x)$ is surjective. However, note that $f(1/p) =0$ and so it is not injective. Hence, by definability of $f(x)$, $(\Z_{p^{\infty}},+)$ is not pseudofinite.
\eej
\bl\label{SumTwoNoPseudofiniteGroups} For any prime $p$, the abelian group $\Z(p^{\infty})\oplus\Z_{(p)}$ is pseudofinite.
\el

\bdem Let $p$ be prime and $\U$ a non principal ultrafilter on the natural numbers, we will prove the following elementary equivalence:
\begin{align*}
    \mathbb{Z}(p^{\infty})\oplus\mathbb{Z}_{(p)}\equiv\prod\limits_{\mathcal{U}} \mathbb{Z}(p^{k}) ,
\end{align*}

where the groups $\Z(p^{k})$ represent the cyclic groups of order $p^k$ and the index set of the ultraproduct varies over $k$. First we will compute the Szmielew invariants of the direct sum in the left hand side for any prime $q$ and $n>0$ in Table \ref{InvariantesSumaDeLosNoPseudofinitos}.

\begin{table}[h!]
\centering
\begin{tabular}{|c|c|c|c|}
\hline
\textbf{Invariant} & $\mathbb{Z}(p^{\infty})$                                                               & $\mathbb{Z}_{(p)}$                                                           & $\mathbb{Z}(p^{\infty})\oplus\mathbb{Z}_{(p)}$                                                             \\ \hline
$\operatorname{U}(q,n;-)$          & 1                                                                         & 1                                                                         & 1                                                                         \\ \hline
$\operatorname{D}(q,n;-)$          & \begin{tabular}[c]{@{}c@{}}$p$ if $(p=q)$\\ $1$ otherwise\end{tabular} & 1                                                                         & \begin{tabular}[c]{@{}c@{}}$p$ if $(p=q)$\\ $1$ otherwise\end{tabular} \\ \hline
$\operatorname{Tf}(q,n;-)$         & 1                                                                         & \begin{tabular}[c]{@{}c@{}}$p$ if $(p=q)$\\ $1$ otherwise\end{tabular} & \begin{tabular}[c]{@{}c@{}}$p$ if $(p=q)$\\ $1$ otherwise\end{tabular} \\ \hline
$\operatorname{Exp}(q,n;-)$        & $\infty$                                                                  & $\infty$                                                                  & $\infty$                                                                  \\ \hline
\end{tabular}
\caption{Table with the Szmielew invariants of the sum.}
\label{InvariantesSumaDeLosNoPseudofinitos}
\end{table}

Now we will show that the Szmielew invariants of the ultraproduct $\prod\limits_{\mathcal{U}} \mathbb{Z}(p^{k})$ are equal to the ones in the direct sum and by Theorem \ref{SzmielewTheorem} we will conclude that these are elementary equivalent.
\bitem
\item For the \textbf{Ulm invariant} we have the following two cases:\\

\underline{\textbf{Case 1: $q=p$}} Then for all $k\neq n+1$, we get that $\mathbb{Z}(p^{k})\models \sigma_{q,n,1}^{\operatorname{U}}$, where the sentence $\sigma_{q,n,1}^{\operatorname{U}}$ characterizes that $\operatorname{U}(p,n; \mathbb{Z}(p^{k}))=1$. Then the following index set $\{k\in\mathbb{N}:\mathbb{Z}(p^{k})\models \sigma_{q,n,1}^{\operatorname{U}}\}$ is in the ultrafilter $\U$ and by \Los's theorem, $\sigma_{q,n,1}^{\operatorname{U}}$ is true in the ultraproduct.\\

\underline{\textbf{Case 2: $q\neq p$}} For all natural $k$, $\mathbb{Z}(p^{k})\models \sigma_{q,n,1}^{\operatorname{U}}$ and so the sentence is true in the ultraproduct.\\

\noindent By the previous cases, we get for all prime $q$ that
\begin{align*}
    \prod\limits_{\mathcal{U}} \mathbb{Z}(p^{k})\models \sigma_{q,n,1}^{\operatorname{U}}. 
\end{align*}
So the Ulm invariant of the ultraproduct is equal to 1 for any prime $q$.

\item We compute the \textbf{Divisible invariant} in the following two cases:\\

\underline{\textbf{Case 1: $q=p$}} Then for all $k> n$, we get that $\mathbb{Z}(p^{k})\models \sigma_{q,n,p}^{\operatorname{D}}$, where the sentence $\sigma_{q,n,p}^{\operatorname{D}}$ characterizes that $\operatorname{D}(p,n; \mathbb{Z}(p^{k}))=p$. Observe that the following index set $\{k\in\mathbb{N}:\mathbb{Z}(p^{k})\models \sigma_{q,n,p}^{\operatorname{D}}\}$ is in the ultrafilter $\U$ and by \Los's theorem, $\sigma_{q,n,p}^{\operatorname{D}}$ is true in the ultraproduct.\\

\underline{\textbf{Case 2: $q\neq p$}} For all natural $k$, $\mathbb{Z}(p^{k})\models \sigma_{q,n,1}^{\operatorname{D}}$ and so the sentence is true in the ultraproduct.\\

\noindent By the previous cases, we conclude that for $q$ prime:
\begin{align*}
    &\text{If $q\neq p$, then }\prod\limits_{\mathcal{U}} \mathbb{Z}(p^{k})\models \sigma_{q,n,1}^{\operatorname{D}}.\hspace{0.5cm}\text{Otherwise }\prod\limits_{\mathcal{U}} \mathbb{Z}(p^{k})\models \sigma_{p,n,p}^{\operatorname{D}}. 
\end{align*}
Finishing the computation of the Divisible invariant of the ultraproduct. Furthermore, the Torsion-free invariant of the ultraproduct is exactly the same since both of these invariants behave similarly for finite cyclic groups.

\item Finally we show that the \textbf{Exponential invariant} is infinite in the ultraproduct by considering the following cases:\\

\underline{\textbf{Case 1: $q=p$}}. Let $\ell$ be the natural number such that for all $\ell'>\ell$ we get that $p^{\ell'}>m$, and define $k_0:=n+\ell+1$. Note that for all $k\geq k_0$ we get that $\mathbb{Z}(p^{k})\not\models \sigma_{q,n,m}^{\operatorname{Exp}}$ for any natural $m<\omega$. Then the following index set $\{k\in\mathbb{N}:\mathbb{Z}(p^{k})\nvDash \sigma_{q,n,m}^{\operatorname{Exp}}\}$ is part of the ultrafilter $\U$ and using \Los's Theorem we conclude that the sentence $\sigma_{q,n,m}^{\operatorname{Exp}}$ is not true in the ultraproduct.\\

\underline{\textbf{Case 2: $q\neq p$}}. By an equivalent argument as before, we obtain that the sentence $\sigma_{q,n,m}^{\operatorname{Exp}}$ is not true in the ultraproduct as well.\\

\noindent By the previous cases, we conclude that for $q$ prime:
\begin{align*}
    \prod\limits_{\mathcal{U}} \mathbb{Z}(p^{k})\not\models \sigma_{q,n,m}^{\operatorname{Exp}}. 
\end{align*}
for any finite natural $m<\omega$. Thus, the Exponential invariant is infinite.
\eitem
As the Szmielew invariants of the abelian group $\mathbb{Z}(p^{\infty})\oplus\mathbb{Z}_{(p)}$ and $\prod\limits_{\mathcal{U}} \mathbb{Z}(p^{k})$ are the same, by Theorem \ref{SzmielewTheorem} these are elementary equivalent and so $\mathbb{Z}(p^{\infty})\oplus\mathbb{Z}_{(p)}$ is pseudofinite.
\edem
\bc For all $\kappa\leq\omega$ and $p$ prime, the abelian group:
\begin{align*}
    \left(\mathbb{Z}(p^{\infty})\oplus\mathbb{Z}_{(p)}\right)^{\kappa},
\end{align*}
is pseudofinite.
\ec
\bdem By Corollary \ref{CountableSumPseudofinite} and Lemma \ref{SumTwoNoPseudofiniteGroups} we conclude that the group is pseudofinite.
\edem
Having established which combinations of building blocks are automatically pseudofinite, we now turn the question around: given an arbitrary Szmielew group, which choices of its parameters $(\kappa_{p,n-1},\lambda_p,\mu_p,\nu)$ make it pseudofinite? The answer comes from comparing two specific Szmielew invariants which are the \textit{Divisible invariant} and the \textit{Torsion free invariant}, since, as we will show next, these two invariants coincide on every finite abelian group, and pseudofinite groups must inherit this property.

\bt\label{DivisibleTorsionInvariansFiniteABGroups} Every finite abelian group satisfies that for all prime $p$ and $n>0$, the \textit{Divisible} and \textit{Torsion free} invariants are equal.
\et
\bdem Let $G$ be a finite abelian group. By the fundamental theorem of finite abelian groups, G is isomorphic to a direct sum of finite cyclic abelian groups of the form:
\begin{align*}
    G\cong \mathbb{Z}(p_{1}^{k_{1}})\oplus\mathbb{Z}(p_{2}^{k_{2}})\oplus\dots\oplus\mathbb{Z}(p_{n}^{k_{n}})
\end{align*}
Where each $\mathbb{Z}(p_{i}^{k_{i}})$ is the cyclic group of order $p_{i}^{k_{i}}$ with $p_i$ prime. Note that for any prime $q$ and $n>0$, $\operatorname{D}(q,n;\mathbb{Z}(p_{i}^{k_{i}})) = \operatorname{Tf}(q,n; \mathbb{Z}(p_{i}^{k_{i}}))$. Now we compute the following invariant of $G$ using the direct sum as follows:
\begin{align*}
    \operatorname{D}(q,n;G)&=\operatorname{D}\left(q,n; \mathbb{Z}(p_{1}^{k_{1}})\oplus\mathbb{Z}(p_{2}^{k_{2}})\oplus\dots\oplus\mathbb{Z}(p_{n}^{k_{n}})\right)\\
    &=\text{sup}\left\{\prod\limits_{j\in J}\operatorname{D}(q,n;\mathbb{Z}(p_{i}^{k_{i}})) : J\subseteq^{fin}\{1,2,\dots,n\}\right\}\\
    &=\text{sup}\left\{\prod\limits_{j\in J}\operatorname{Tf}(q,n;\mathbb{Z}(p_{i}^{k_{i}})) : J\subseteq^{fin}\{1,2,\dots,n\}\right\}\\
    &=\operatorname{Tf}\left(q,n; \mathbb{Z}(p_{1}^{k_{1}})\oplus\mathbb{Z}(p_{2}^{k_{2}})\oplus\dots\oplus\mathbb{Z}(p_{n}^{k_{n}})\right) = \operatorname{Tf}(q,n;G).
\end{align*}
Therefore, for all prime $q$, $n>0$ and finite abelian group $G$, $\operatorname{D}(q,n;G) = \operatorname{Tf}(q,n;G)$. 
\edem
\bo\label{ObsFINTheory} Let $\Gamma$ be the following set of sentences
\begin{align*}
    \Gamma = \bigcup\limits_{p \text{ prime}}\left(\{\tau_{n,k}^{p}: n,k\in\mathbb{N}\}\right),
\end{align*}
where $\tau_{n,k}^{p}$ is defined as:
\begin{align*}
    \tau_{n,k}^{p}: \sigma_{p,n,k}^{\operatorname{D}}\Longleftrightarrow\sigma_{p,n,k}^{\operatorname{Tf}}.
\end{align*}
for $n>0$, $k\in\mathbb{N}$ and $p$ prime. By Theorem \ref{DivisibleTorsionInvariansFiniteABGroups}, $\Gamma$ is part of the common theory of finite abelian groups, or $\text{FIN}_{\Le}$.
\eo
\bl\label{DivisibleTorsionInvariantsOfASzmielewGroup} Let $G=\bigoplus\limits_{p \text{ prime}} \left[ \bigoplus\limits_{n>0} \mathbb{Z}(p^{n})^{\kappa_{(p,n-1)}} \oplus \mathbb{Z}(p^{\infty})^{\lambda_{p}} \oplus \mathbb{Z}^{\mu_p}_{(p)} \right] \oplus \mathbb{Q}^{\nu}$ be a Szmielew group, $p$ prime and $n>0$. Then:
\begin{align*}
    &\operatorname{D}(p,n,G)=p^{\sum\limits_{n<k}\kappa_{p,k-1}}\cdot p^{\lambda_{p}},\hspace{0.3cm}\operatorname{Tf}(p,n,G)=p^{\sum\limits_{n<k}\kappa_{p,k-1}}\cdot p^{\mu_{p}}.
\end{align*} 
\el
\bdem Let $G$ be a Szmielew group, $p$ prime and $n>0$ as in the statement. Let $A_p=\bigoplus\limits_{k>0} \mathbb{Z}(p^{k})^{\kappa_{(p,k-1)}} \oplus \mathbb{Z}(p^{\infty})^{\lambda_p} \oplus \mathbb{Z}_{(p)}^{\mu_{p}}$, so $G=\mathbb{Q}^{\nu}\oplus\bigoplus\limits_{p \text{ prime}} \left[ A_{p} \right]$. Since $\operatorname{D}(p,n,\mathbb{Q}^{\nu})=\operatorname{Tf}(p,n,\mathbb{Q}^{\nu})=1$, and for all prime $q$ different from $p$, $\operatorname{D}(p,n,A_{q})=\operatorname{Tf}(p,n,A_{q})=1$, the behavior of the invariant for the prime $p$ depends only on the factor $A_p$.\\
First we start computing the divisible invariant to obtain:
\begin{align*}
    \operatorname{D}(p,n,A_p)&=\operatorname{D}\left(p,n,\bigoplus\limits_{k>0} \mathbb{Z}(p^{k})^{\kappa_{(p,k-1)}} \oplus \mathbb{Z}(p^{\infty})^{\lambda_p} \oplus \mathbb{Z}_{(p)}^{\mu_{p}}\right)\\
    % &=\operatorname{D}\left(p,n,\bigoplus\limits_{k>0} \mathbb{Z}(p^{k})^{\kappa_{(p,k-1)}}\right)\cdot \operatorname{D}\left(p,n,\mathbb{Z}(p^{\infty})^{\lambda_p}\right) \cdot \operatorname{D}\left(p,n,\mathbb{Z}_{(p)}^{\mu_{p}}\right)\\
    &=\operatorname{D}\left(p,n,\bigoplus\limits_{k>0} \mathbb{Z}(p^{k})^{\kappa_{(p,k-1)}}\right)\cdot \operatorname{D}\left(p,n,\mathbb{Z}(p^{\infty})^{\lambda_p}\right).
\end{align*}
Now we compute the invariant of each one of these factors, first note that $\operatorname{D}\left(p,n,\mathbb{Z}(p^{\infty})\right)=p$, and so $\operatorname{D}\left(p,n,\mathbb{Z}(p^{\infty})^{\lambda_p}\right)=p^{\lambda_{p}}$. On the other hand, $\operatorname{D}\left(p,n,\mathbb{Z}(p^{k})\right)=p$ if $n<k$, otherwise it is equal to $1$. However, the invariant of a direct sum takes the supremum over all finite products of invariants, so the invariant of this factor is given by $\operatorname{D}\left(p,n,\bigoplus\limits_{k>0} \mathbb{Z}(p^{k})^{\kappa_{(p,k-1)}}\right)=p^{\sum\limits_{n<k}\kappa(p,k-1)}$. Therefore we can conclude that the invariant of the factor $A_p$ is given by:
\begin{align*}
    \operatorname{D}(p,n,A_p)=p^{\sum\limits_{n<k}\kappa_{p,k-1}}\cdot p^{\lambda_{p}}.
\end{align*}
Furthermore, note that the computation of the torsion free invariant can be reduced to:
\begin{align*}
    \operatorname{Tf}(p,n,A_p)&=\operatorname{Tf}\left(p,n,\bigoplus\limits_{k>0} \mathbb{Z}(p^{k})^{\kappa_{(p,k-1)}} \oplus \mathbb{Z}(p^{\infty})^{\lambda_p} \oplus \mathbb{Z}_{(p)}^{\mu_{p}}\right)\\
    % &=\operatorname{Tf}\left(p,n,\bigoplus\limits_{k>0} \mathbb{Z}(p^{k})^{\kappa_{(p,k-1)}}\right)\cdot \operatorname{Tf}\left(p,n,\mathbb{Z}(p^{\infty})^{\lambda_p}\right) \cdot \operatorname{Tf}\left(p,n,\mathbb{Z}_{(p)}^{\mu_{p}}\right)\\
    &=\operatorname{Tf}\left(p,n,\bigoplus\limits_{k>0} \mathbb{Z}(p^{k})^{\kappa_{(p,k-1)}}\right)\cdot \operatorname{Tf}\left(p,n,\mathbb{Z}_{(p)}^{\mu_{p}}\right).
\end{align*}
Note that $\operatorname{Tf}\left(p,n,\mathbb{Z}_{(p)}^{\mu_{p}}\right)=p^{\mu_{p}}$, and by observation \ref{ObsFINTheory},
\begin{align*}
    \operatorname{Tf}\left(p,n,\bigoplus\limits_{k>0} \mathbb{Z}(p^{k})^{\kappa_{(p,k-1)}}\right)=p^{\sum\limits_{n<k}\kappa(p,k-1)}.
\end{align*}
So we obtain the following formula:
\begin{align*}
    \operatorname{Tf}(p,n,A_p)=p^{\sum\limits_{n<k}\kappa_{p,k-1}}\cdot p^{\mu_{p}}.
\end{align*}
\edem 
\bd[Divergence of exponent series] For the exponents $\kappa_{p,k-1}$ which appear as powers of the finite cyclic groups in a Szmielew group, we say that the series $\sum\limits_{k>0}\kappa_{p,k-1}$ \textbf{diverges} if there are infinitely many natural numbers for which $\kappa_{p,k-1}\neq0$. Otherwise we say that the series \textbf{converges}.
\ed
\bl\label{PseudofinitePowerDivergence} Let $p$ be a prime number and consider the group $A_p=\bigoplus\limits_{k>0} \mathbb{Z}(p^{k})^{\kappa_{(p,k-1)}} \oplus \mathbb{Z}(p^{\infty})^{\lambda_p} \oplus \mathbb{Z}_{(p)}^{\mu_{p}}$. If $\sum\limits_{k>0}\kappa_{p,k-1}$ diverges or $\lambda_{p}=\mu_{p}$, then $A_p$ is pseudofinite.
\el
\bdem First assume that $\sum\limits_{k>0}\kappa_{p,k-1}$ diverges and let $q$ be a prime different from $p$, $n>0$. Using the results in the proof of Lemma \ref{DivisibleTorsionInvariantsOfASzmielewGroup}, we get the following results for the divisible invariant:
\begin{align*}
    \operatorname{D}(p,n,A_p)=p^{\sum\limits_{n<k}\kappa_{p,k-1}}\cdot p^{\lambda_{p}}=\infty,\hspace{0.5cm}\operatorname{D}(q,n,A_p)=1.
\end{align*}
We can conclude a similar result for the torsion free invariant as follows:
\begin{align*}
    \operatorname{Tf}(p,n,A_p)=p^{\sum\limits_{n<k}\kappa_{p,k-1}}\cdot p^{\mu_{p}}=\infty,\hspace{0.5cm}\operatorname{Tf}(q,n,A_p)=1.
\end{align*}
On the other hand, note that the Ulm invariant only depends on the sum of the finite cyclic groups of the form:
\begin{align*}
    \operatorname{U}(p,n,A_p)=\operatorname{U}\left(p,n,\bigoplus\limits_{k>0} \mathbb{Z}(p^{k})^{\kappa_{(p,k-1)}} \right)=p^{\kappa_{p,n-2}},\hspace{0.5cm}\operatorname{U}(q,n,A_p)=1.
\end{align*}
Furthermore, as $A_p$ has two infinite summands, we get that $\operatorname{Exp}(p,n,A_p)=\operatorname{Exp}(q,n,A_p)=\infty$.\\

Moreover, we will show that $A_p$ has the same Szmielew invariants as $H=\bigoplus\limits_{k>0} \mathbb{Z}(p^{k})^{\kappa_{(p,k-1)}}$. Note that we just have to show that the exponential invariant of $H$ diverges, which we compute as follows:
\begin{align*}
    \operatorname{Exp}\left(p,n,H\right) =\operatorname{Exp}\left(p,n,\bigoplus\limits_{k>0} \mathbb{Z}(p^{k})^{\kappa_{(p,k-1)}}\right) &= \prod\limits_{n<k}(p^{k-n})^{\kappa_{p,k-1}}\cdot\prod\limits_{n>k} (p^{k})^{\kappa_{p,k-1}}=\infty.\\
    \operatorname{Exp}\left(p,n,H\right) =\operatorname{Exp}\left(q,n,\bigoplus\limits_{k>0} \mathbb{Z}(p^{k})^{\kappa_{(p,k-1)}}\right) &= \prod\limits_{k>0}(p^{k})^{\kappa_{(p,k-1)}}=\infty.
\end{align*}
By Theorem \ref{SzmielewTheorem}, $A_p\equiv H$ and so it is pseudofinite.\\

Now, if $\lambda_{p}=\mu_{p}$, we get that $A_p=\bigoplus\limits_{k>0} \mathbb{Z}(p^{k})^{\kappa_{(p,k-1)}} \oplus (\mathbb{Z}(p^{\infty})\oplus \mathbb{Z}_{(p)})^{\mu_{p}}$. By Lemma \ref{SumTwoNoPseudofiniteGroups}, $A_p$ is a direct sum of pseudofinite abelian groups, and by Theorem \ref{SumaDirectaDosGAPseudo} we conclude that it is pseudofinite.
\edem
\bt[Classification of pseudofinite abelian groups] \label{Theorem:ClassificationPsfAbelianGroups} Let $G$ be an abelian group. The following are equivalent:
\benum
\item $G$ is pseudofinite.
\item For all prime $p$ and $n>0$, $\operatorname{D}(p,n;G) = \operatorname{Tf}(p,n;G)$.
\item Any Szmielew group $H=\bigoplus\limits_{p \text{ prime}} \left[ \bigoplus\limits_{n>0} \mathbb{Z}(p^{n})^{\kappa_{(p,n-1)}} \oplus \mathbb{Z}(p^{\infty})^{\lambda_{p}} \oplus \mathbb{Z}^{\mu_p}_{(p)} \right] \oplus \mathbb{Q}^{\nu}$ elementarily equivalent to $G$ satisfies for all prime $p$ that either $\lambda_{p} = \mu_{p}$, or $\sum\limits_{k>0}\kappa_{p,k-1}$ diverges.
\eenum
\label{ClassificationTheoremForAbelianGroups}
\et
\bdem For this proof, let $G$ be an abelian group, and $H$ be an Szmielew group such that $G\equiv H$.\\

\noindent \underline{\textbf{$(1)\to (2)$}:} If $G$ is pseudofinite, then it is a model of the common theory of finite abelian groups and by Observation \ref{ObsFINTheory}, $G$ satisfies for any prime $p$ and $n>0$, $\operatorname{D}(p,n;G) = \operatorname{Tf}(p,n;G)$.\\

\noindent \underline{\textbf{$(2)\to (3)$}:} By Lemma \ref{DivisibleTorsionInvariantsOfASzmielewGroup}, the divisible and torsion free invariants can be explicitly computed as:
\begin{align*}
    \operatorname{D}(p,n,G)=p^{\sum\limits_{n<k}\kappa_{p,k-1}}\cdot p^{\lambda_{p}},\hspace{0.5cm}\operatorname{Tf}(p,n,G)=p^{\sum\limits_{n<k}\kappa_{p,k-1}}\cdot p^{\mu_{p}}.
\end{align*} 
If we assume that these invariants are equal for all $p$ prime and $n>0$, then either the sum $\sum\limits_{k>0}\kappa_{p,k-1}$ converges and $\lambda_p=\mu_p$, or the sum $\sum\limits_{k>0}\kappa_{p,k-1}$ diverges.\\

\noindent \underline{\textbf{$(3)\to (1)$}:} If for all prime $p$ that either $\lambda_{p} = \mu_{p}$, or $\sum\limits_{n>0}\kappa_{p,k-1}$ diverges, by Lemma \ref{PseudofinitePowerDivergence} each factor $A_p$ of $H$ is pseudofinite. Then, $H$ is a countable direct sum of pseudofinite factors and by Corollary \ref{CountableSumPseudofinite} we can conclude that $H$ is pseudofinite which implies that $G$ is pseudofinite as well.
\edem

\subsection{Pseudofinite $\omega$-stable abelian groups}
\label{Subsec:PsfWStableAG}
Theorem \ref{ClassificationTheoremForAbelianGroups} characterizes pseudofiniteness of an abelian group in terms of its Szmielew invariants. It is natural to ask how this notion relates to the model-theoretic complexity of the group, in particular to its $\omega$-stability. Recall that Macintyre gave a clean structural classification of $\omega$-stable abelian groups, which we will recall in Theorem \ref{TeoremaMacintyre2.0}, that looks similar to our given results in Theorem \ref{ClassificationTheoremForAbelianGroups}.\\

Hence, we will recall the definitions of $\omega$-stability and show with two explicit examples that pseudofiniteness and $\omega$-stability are independent properties. Later, we will work towards a classification theorem that combines the properties of an abelian group being both pseudofinite and $\omega$-stable (Theorem \ref{ClassificationWStablePsfAbGroups})
\bd[$\omega$-stable theories] A theory $T$ is $\omega$-stable if for every model $\mathcal{M}$ of $T$, if $A$ is a countable subset of $M$, then $|S_1(A)|\leq\omega$. (i.e. the set of 1-types over $A$ is countable) 
\ed
For $\omega$-stable abelian groups, the following chain property for definable subgroups is satisfied.
\bl If $G$ is an $\omega$-stable abelian group, there is no definable collection of subgroups of $G$, $\{H_{i}\}_{i\in\mathbb{N}}$ such that $H_{1}\gneq H_{2}\gneq\cdots H_{i}\gneq H_{i+1}\cdots$. (i.e. there is no infinite descending chain of definable subgroups of $G$)
\label{CadenasDescenInfinitas}
\el
\bd Let $G$ be an abelian group:
\begin{itemize}
    \item $G$ is \textit{divisible} if for all $g\in G$ and $n>1$, there is some $h\in G$ such that $G\models n\cdot h=g$.
    \item $G$ has \textit{bounded exponent} if there is a positive integer $n$ such that $nG = \{0\}$.
\end{itemize}
\ed

The following lemmas represent a characterization of both divisible and bounded exponent abelian groups.
\bl[See Theorem 4 in \cite{kaplansky2018infinite}] An abelian group $D$ is divisible if and only if
\begin{align*}
    D\cong\mathbb{Q}^{\beta}\oplus\bigoplus\limits_{p \text{ prime}}(\mathbb{Z}(p^{\infty})^{d_{p}}).
\end{align*}
for some cardinals $\beta,d_{p}$.
\label{IsomorAbelianosDivisibles}
\el

\bl[See Theorem D in \cite{Macintyre1971}] If $H$ is an abelian group with bounded exponent, then:
\begin{align*}
    H \cong \bigoplus\limits_{p \text{ prime}}\left(\bigoplus\limits_{n>0} \mathbb{Z}(p^{n})^{\gamma_{p,n}}\right).
\end{align*}
where each $\gamma_{p,n}$ is a cardinal and only finitely many exponents $\gamma_{p,n}$ are non-zero.
\label{IsomorAbelExponenteAcotado}
\el 
The following corollary shows that there is a relationship between pseudofiniteness and bounded exponent abelian groups.
\bc Every countable abelian group with bounded exponent is pseudofinite. 
\ec 
With the previous definitions and lemmas, we introduce Macintyre's classification theorem.
\bt[Macintyre's Theorem (see Theorem 2.1 in \cite{Macintyre1971})] Let $G$ be an abelian group. $\operatorname{Th}(G)$ is $\omega$-stable if and only if $G$ is isomorphic to $D\oplus H$ where $D$ is a divisible abelian group and $H$ is an abelian group with bounded exponent. 
\label{TeoremaMacintyre2.0}
\et
Notice that with Theorem \ref{TeoremaMacintyre2.0} we have a way to characterize $\omega$-stable abelian groups, and using our previous results from Theorem \ref{ClassificationTheoremForAbelianGroups} we can classify pseudofiniteness of abelian groups. Hence, the remaining part of this section focus on the possibility of combining both theorems to relate pseudofiniteness with $\omega$-stability in abelian groups. However, we shall observe that there is no immediate relationship between both concepts in our setup. For the following two examples, let $p$ be a prime number.
\bej Consider the abelian group $H=\mathbb{Z}(p^{\infty})\oplus \mathbb{Z}(p)$ which is $\omega$-stable by Theorem \ref{TeoremaMacintyre2.0} as $\mathbb{Z}(p^{\infty})$ is divisible and $\mathbb{Z}(p)$ has bounded exponent. (In particular, the bound is $p$) However, $H$ is not pseudofinite as the following invariants do not agree for the same prime $p$ and $n=2$:
\begin{align*}
    \operatorname{D}(p,2,H) &= \operatorname{D}(p,2,\mathbb{Z}(p^{\infty}))\cdot\operatorname{D}(p,2,\mathbb{Z}(p)) = \operatorname{D}(p,2,\mathbb{Z}(p^{\infty})) = p.\\
    \operatorname{Tf}(p,2,H) &= \operatorname{Tf}(p,2,\mathbb{Z}(p^{\infty}))\cdot\operatorname{Tf}(p,2,\mathbb{Z}(p)) = 1.
\end{align*}
Hence, $\operatorname{D}(p,2,H)\neq\operatorname{Tf}(p,2,H)$, and by Theorem \ref{ClassificationTheoremForAbelianGroups}, $H$ is not pseudofinite.
\eej
\bej Consider the abelian group $G=\mathbb{Z}(p^{\infty})\oplus\mathbb{Z}_{(p)}$ which is pseudofinite as shown in Lemma \ref{SumTwoNoPseudofiniteGroups}. However, for all $i\geq1$ we can consider the formula $\theta_{i}(x):=p^{i}|x$. We will show that this formula defines an infinite descending chain of subgroups of $G$. First, note that the subgroups defined by each $\theta_i$ are the following:
\begin{align*}
    \theta_{i}(\mathbb{Z}(p^{\infty})\oplus\mathbb{Z}_{(p)})=\{x\in\mathbb{Z}(p^{\infty})\oplus\mathbb{Z}_{(p)}: \mathbb{Z}(p^{\infty})\oplus\mathbb{Z}_{(p)}\models\theta_{i}(x)\} = \mathbb{Z}(p^{\infty})\oplus p^{i}\mathbb{Z}_{(p)},
\end{align*}
Since $\mathbb{Z}(p^{\infty})\oplus p^{i}\mathbb{Z}_{(p)}\geq\mathbb{Z}(p^{\infty})\oplus p^{i+1}\mathbb{Z}_{(p)}$ for all $i\geq1$, we obtain an infinite descending chain of definable subsets of $G$ and by Lemma \ref{CadenasDescenInfinitas}, we conclude that $G$ is not $\omega$-stable.
\eej
The previous examples show that there is no direct relationship between these two concepts. However, we can formulate our classification theorem for pseudofinite $\omega$-stable abelian groups as follows.
\bt[Classification of pseudofinite $\omega$-stable abelian groups] \label{ClassificationWStablePsfAbGroups}An abelian group $G$ is pseudofinite and $\omega$-stable if and only if $G$ is isomorphic to an abelian group of the form:
\begin{align*}
    \mathbb{Q}^{\beta}\oplus\bigoplus\limits_{i=1}^{K}\mathbb{Z}(p_{i}^{n_{i}})^{\gamma_{p_{i},n_{i}}},
\end{align*}
for $K$ natural, $p_{i}$ prime, $n_{i}>0$ natural, and $\beta$, $\gamma_{p,n}$ cardinals less than or equal to $\omega$.
\et
\bdem Let $G$ be an abelian group. Assume that $G \cong \mathbb{Q}^{\beta}\oplus\bigoplus\limits_{i=1}^{K}\mathbb{Z}(p_{i}^{n_{i}})^{\gamma_{p_{i},n_{i}}}$ where $\beta$, $\gamma_{p,n}$ are cardinals less than or equal to $\omega$. By Lemma \ref{IsomorAbelExponenteAcotado}, the direct sum of finitely many cyclic groups has bounded exponent. So, by Theorem \ref{TeoremaMacintyre2.0}, $G$ is $\omega$-stable and, since each summand is pseudofinite, using Corollary \ref{CountableSumPseudofinite} we conclude that $G$ is pseudofinite.\\

On the other hand, assume that $G$ is $\omega$-stable and pseudofinite. As $G$ is $\omega$-stable, Theorem \ref{TeoremaMacintyre2.0} implies that $G\cong D\oplus H$ where $D$ is a divisible abelian group and $H$ is an abelian group with bounded exponent. Then, using Lemma \ref{IsomorAbelianosDivisibles}, we get that $D\cong\mathbb{Q}^{\beta}\oplus\bigoplus\limits_{p}(\mathbb{Z}(p^{\infty})^{d_{p}})$, and by Lemma \ref{IsomorAbelExponenteAcotado}, we can write $H \cong \bigoplus\limits_{p  \text{ prime}}\left(\bigoplus\limits_{n>0} \mathbb{Z}(p^{n})^{\gamma_{p,n}}\right)$ 
where each $\gamma_{p,n}$ is a cardinal and only finitely many exponents $\gamma_{p,n}$ are non-zero. This gives us the following isomorphism for $G$.
\begin{align*}
    G\cong D\oplus H \cong \mathbb{Q}^{\beta}\oplus\bigoplus\limits_{p}(\mathbb{Z}(p^{\infty})^{d_{p}}) \oplus \bigoplus\limits_{p  \text{ prime}}\left(\bigoplus\limits_{n>0} \mathbb{Z}(p^{n})^{\gamma_{p,n}}\right)
\end{align*}
However, since $G$ is pseudofinite and there is no summand of the form $\mathbb{Z}_{(p)}$ in $H$, by the Classification of pseudofinite abelian groups (Theorem \ref{ClassificationTheoremForAbelianGroups}) we get that $d_{p}=0$ for all prime $p$. We conclude that $G$ is isomorphic to a group of the following form:
\begin{align*}
    G\cong \mathbb{Q}^{\beta} \oplus \bigoplus\limits_{p  \text{ prime}}\left(\bigoplus\limits_{n>0} \mathbb{Z}(p^{n})^{\gamma_{p,n}}\right)
\end{align*}
which matches the one in our theorem statement since only finitely many exponents $\gamma_{p,n}$ are non-zero in this construction.
\edem

\section{Classes of finite abelian groups and measurability}
The guiding question is which classes of finite abelian groups are tame enough for the cardinalities of their definable sets to be controlled uniformly, in the sense of multidimensional exact classes of Anscombe, Macpherson, Steinhorn and Wolf in \cite{AMSW}, and of the $N$-dimensional asymptotic classes defined by Macpherson and Steinhorn, and of Elwes in \cite{Macpherson-Steinhorn,Elwes-AsymptoticClasses}. In Subsection \ref{SS:MEC} we show that the class of finite abelian groups forms a m.e.c (see Proposition \ref{Prop:MecAG}). This result is already in the literature (see Theorem 4.2.3 in \cite{AMSW}), but our proof, based on the Smith normal form, describes more
explicitly the possible cardinalities of definable sets.\\

In Subsection \ref{SS:MS-measurability} we use Roger's characterization of superstable abelian groups (see Theorem \ref{thm:Rogers}) to show that pseudofiniteness does not imply MS-measurability, which contradicts a claim of Kestner in \cite[Remark 4.14]{Kestner2014}. Finally, in Subsection \ref{SS:AsymptoticClasses}, we address a possible characterization of $N$-dimensional asymptotic classes of finite abelian groups. The recurring obstruction is that ultraproducts of $N$-dimensional asymptotic classes are supersimple of finite $SU$-rank, which for abelian groups forces superstability. We introduce the notions of width, depth and entangled classes to organize the examples, and we close with a conjecture on the general picture (Conjecture \ref{Conj:NDimAsympt}).
\subsection{Alternative proof on the multidimensional exactness of the class of finite abelian groups}\label{SS:MEC}
The definition of multidimensional exact classes appear in \cite{Wolf-MEC} and \cite{AMSW}. It is a way to encompass a rich family of classes of finite structures with strong conditions on the possible cardinalities of definable sets, without restricting too much the possible values for dimensions and measures, and allowing ultraproducts that may have infinite SU-rank. We now collect here the definition and the main results regarding exact classes of finite structures.
\bd [$\mathcal{R}$-mec, alternative definition] \label{Def:MEC} Let $\mathcal{C}$ be a class of finite $\Le$-structures, and let $\mathcal{R}$ be a set of functions from $\mathcal{C}$ to $\mathbb{N}$. We say that $\mathcal{C}$ is a \emph{multidimensional exact class for $\mathcal{R}$ in $\Le$} (or \emph{$\mathcal{R}$-mec in $\mathcal{L}$} for short) if for every parameter-free $\Le$-formula $\varphi(\ov{x},\ov{y})$ with $m=|\ov{y}|$ and $n=|\ov{x}|$ there exist finitely many functions $h_1,\ldots,h_k\in \mathcal{R}$ and formulas $\psi_1(\ov{y}),\ldots,\psi_k(\ov{y})$ such that the following hold:
\benum
\item For every $\mathcal{M}\in \mathcal{C}$, the formulas $\psi_1(\ov{y}),\ldots,\psi_k(\ov{y})$ provide a definable partition of $\mathcal{M}^m$.
\item For every structure $\mathcal{M}$ in $\mathcal{C}$ and every tuple $\ov{a}$ from $\mathcal{M}^m$, we have \[|\varphi(\mathcal{M};\ov{a})|=h_i(\mathcal{M}) \text{\ if and only if\ }\mathcal{M}\models \psi_i(\ov{a}).\]
\eenum
Also, we say that $\mathcal{C}$ is an \emph{$R$-polynomial exact class} if the set of functions $\mathcal{R}$ is the collection of polynomials with coefficients in some ring $R$, possibly with several variables. In the case where $R=\mathbb{Q}$, we simply call $\mathcal{C}$ a polynomial exact class.
\ed
\brmk In the previous definition, the fact that there are finitely many functions is often referred to as the \emph{finiteness condition}, while the existence of the formulas $\psi_i(\ov{y})$ describing the cardinalities of definable sets is referred to as the \emph{definability condition}. 
\ermk

The fact that these are called multidimensional \emph{exact} classes refers to the fact that the cardinality can be exactly described by one of the functions $h_i$, as opposed to the multidimensional \emph{asymptotic} classes where the cardinality of definable sets is described asymptotically by one of the functions $h_i$, with a small error that goes to zero as we consider arbitrary large finite structures. In this work we will restrict ourselves to examples of multidimensional exact classes. For the more general concept of m.a.c.'s and their systematic study we refer the reader to \cite{AMSW}.\\

By the results of A. Pillay in \cite{Pillay-stronglyminimalpsf} we know that if $\mathcal{C}$ is a class of finite $\Le$-structures and all infinite ultraproducts of structures in $\mathcal{C}$ are strongly minimal, then $\mathcal{C}$ is a polynomial exact class and the functions $h$ correspond to polynomials in one variable and integer coefficients. This includes for example the class $\mathcal{C}=\{(\mathbb{Z}/p\mathbb{Z},+):\text{$p$ prime}\}$, as all infinite ultraproducts are models of the theory of torsion-free divisible abelian groups. Pillay's results were generalized by A. van Abel in \cite{vanAbel-counting} to classes whose infinite ultraproducts are uncountably categorical, and in this case the polynomials may have rational coefficients. \\

Some other examples of m.e.c include: 
\bitem
\item As presented in \cite{AMSW}, let $\mathcal{C}_d$ be the class of all structures $(M,E)$ where $E$ is an equivalence relation on $M$ with at most $d$ different sizes of equivalence classes. Then $\mathcal{C}_d$ is a polynomial m.e.c., after expansion by unary predicates picking out the union of the equivalence classes of given size.
\item If $\mathcal{C}_d=\{G_n:n<\omega\}$ is a class of finite $d$-graphs such that $\operatorname{girth}(G_n)\to \infty$ then all infinite ultraproducts of graphs in $\mathcal{C}$ are acyclic $d$-regular graphs. Hence, these ultraproducts are all strongly minimal, so $\mathcal{C}$ is a polynomial exact class.
\item Suppose now that $\mathcal{C}=\{G_n:n<\omega\}$ is a class of regular graphs such that both the degree of regularity $\operatorname{deg}(G_n)$ and the girth $\operatorname{girth}(G_n)$ tend to infinity. It is shown in \cite{Garcia-Robles} that any such class is a multidimensional exact class, and the sizes of definable sets have the form $p(\operatorname{deg}(G_n),\operatorname{girth}(G_n))$, where $p(T_1,T_2)\in \mathbb{Z}[T_1,T_2]$. Infinite ultraproducts of structures in these classes of graphs satisfy the theory of the everywhere infinite forest, so they are $\omega$-stable of Morley rank $\omega$.
\item By Example 2.4.3 of \cite{Wolf-MEC} (based on Theorem 3.14 of \cite{Macpherson-Steinhorn}, the class $\mathcal{C}$ of finite cyclic groups is a m.e.c. 
\item By Proposition 4.4.2 of \cite{Garcia-Macpherson-Steinhorn}, the class $\mathcal{C}$ of all finite homocyclic groups, which are groups of the form $(\Z/p^n\Z)^m$, where $p$ is prime and $m,n$ are positive integers, is a multidimensional exact class whose ultraproducts are strictly stable.

\eitem

A priori, worrying about formulas with more than one variable from the structure and allowing multiple parameters would complicate the computations. In the case of abelian groups, the result from Lemma \ref{PruferPPF} would simplify to solutions of linear equations in a single variable instead of matrix systems of equations with multiple variables. D. Wolf in \cite{Wolf-MEC} showed that proving the existence of measuring functions for single variable formulas was enough to conclude the property of being a m.e.c. as presented in the following fact.

\bh [Projection Lemma, Lemma 2.3.1 in \cite{Wolf-MEC}] \label{Fact:Projection-Lemma} Let $\mathcal{C}$ be a class of $\Le$-structures such that the definition of an $R$-mec holds in $\mathcal{C}$ for all formulas $\varphi(x,\ov{y})$ with a single object variable $x$. Then $\mathcal{C}$ is an $R'$-mec in $\Le$, where $R'$ is generated under addition and multiplication by the functions in $R$.
\eh 

The following results show that, in fact, the class of all finite abelian groups corresponds to a multidimensional exact class. A proof of this result appears in \cite[Theorem 4.2.3]{AMSW}, as a corollary of the more general fact that for every ring $R$, the class of finite $R$-modules is a multidimensional exact class. Even though their proof is correct, it does not provide much information about the possible sizes of definable sets that one may encounter in classes of abelian groups. We provide here an alternative proof, with a longer description of the possible cardinalities of definable sets in abelian groups, and this description will later be used to show that certain classes of finite abelian groups form in fact $N$-dimensional asymptotic classes.

\bl \label{Lemma:KernelSNF} Let $A$ be finite abelian group, and let $T_A:A^n\to A^m$ be the linear map defined by multiplication by an $m\times n$ matrix $T$ with integer coefficients. Then, there are integers $c_1,\ldots,c_\nu$ (depending only on the entries of $T$) such that
\[|\operatorname{ker}(T_A)|=|A|^{n-\nu}\prod_{i=1}^\nu |A/c_iA|.\]
\el
\bdem Using Smith Normal Form, we can write $UTV=D$ for some invertible matrices $U,V$ and a left diagonal matrix $D$ as in Fact \ref{Fact:SmithNormalForm}. Since $U$ and $V$ are invertible, we have

\begin{align*}
|\operatorname{ker}(T_A)|&=|\{\bm{x}\in A^n:T\cdot \bm{x}=\bm{0}\}|=|\{\bm{x}\in A^n:UT\cdot \bm{x}=\bm{0}\}|=|\{\bm{y}\in A^n: UT(V\bm{y})=\bm{0}\}|\\
&=|\{\bm{y}\in A^n: D\cdot \bm{y}=\bm{0}\}|=|\operatorname{ker}(D_A)|.
\end{align*}

On the other hand, since $c_1,\cdots,c_\nu\geq 1$, we have
\begin{align*}
\operatorname{ker}(D_A)&=\left\{\bm{y}=\begin{bmatrix}y_1\\\vdots\\y_n\end{bmatrix}:\begin{bmatrix}
  \begin{matrix}
  c_1 & 0 & \cdots & 0 \\
  0 & c_2 & \cdots & 0 \\
  \vdots & \vdots & \ddots & \vdots\\
  0 & 0 &\cdots & c_{\nu}
  \end{matrix}
  & \rvline & & \bigzero & \\
\hline
& \rvline &\\
  \bigzero & \rvline & & \bigzero
\end{bmatrix}\cdot \begin{bmatrix}y_1\\\vdots\\y_n\end{bmatrix} = \begin{bmatrix}0\\\vdots\\0\end{bmatrix} \right\}\\
&=\{(y_1,y_2,\ldots,y_\nu,\ldots,y_n): c_1\cdot y_1=0, c_2\cdot y_2=0,\ldots,c_\nu\cdot y_\nu=0\}\\
&=\operatorname{ker}(\lambda_{c_1})\times \cdots \times \operatorname{ker}(\lambda_{c_\nu})\times A^{n-\nu}
\end{align*}
Therefore, we have
\begin{align*}
|\operatorname{ker}(T_A)|&=|\operatorname{ker}(D_A)|=|\operatorname{ker}(\lambda_{c_1})\times \cdots \times \operatorname{ker}(\lambda_{c_\nu})\times A^{n-\nu}|\\
&=|A^{n-\nu}|\cdot \prod_{i=1}^\nu |\operatorname{ker}(\lambda_{c_i})|=|A|^{n-\nu}\cdot \prod_{i=1}^\nu \left|A/c_iA\right|.
\end{align*}
\edem

\bl \label{Lemma:sizes-pp-definable} Let $A$ be an abelian group, and consider a formula of the form
\[\theta(x;\ov{y})\,\,:\,\,\bigwedge_{i\in I}[\,k_i\cdot x+t_i(\ov{a})=0\,] \wedge \bigwedge_{j\in J}[\ p_j^{n_j}|(r_j\cdot x+t_j(\ov{y}))\,].\]
\benum
\item[(a)] For every $\ov{a}\in A^{|\ov{y}|}$, either $\theta(A;\ov{a})$ is empty or $\theta(x;\ov{a})$ defines a set of the form $b+\theta(A;0)$, where $\theta(A;\bm{0})$ is a subgroup of $A$.
\item[(b)] There are integers $c_1,\ldots,c_\nu,d_1,\ldots,d_\mu$ (depending only on $\{k_i:i\in I\}\cup \{r_j:j\in J\}$) with $\nu\leq |J|+1, \mu\leq |J|, \nu\leq \mu+1$, such that for every finite abelian group $A$ and $\ov{a}\in A^{|\ov{y}|}$, either $\theta(A;\ov{a})$ is empty or
\[|\theta(A;\ov{a})|=|A|^{1+\mu-\nu}\cdot \dfrac{\prod_{\ell=1}^\nu |A/c_\ell A|}{\prod_{\ell=1}^\mu|A/d_\ell A|}\ \ .\]
\eenum
\el
\bdem Notice that $\theta(A;\bm{0})$ corresponds to the set
\[\left\{\alpha\in A:\bigwedge_{i\in I}[\,k_i\cdot \alpha=0\,] \wedge \bigwedge_{j\in J}[\,p_j^{n_j}|r_j\cdot \alpha\,]\right\}
=\bigcap_{i\in I}\operatorname{ker}(\lambda_{k_i})\cap \bigcap_{j\in J}\lambda^{-1}_{r_j}(p_j^{n_j}A),\]
which is an intersection of subgroups of $A$. Hence, for every $\ov{a}\in A^{|\ov{y}|}$, $\theta(A;\ov{a})$ corresponds to the intersections of cosets of the form $b_i+\operatorname{ker}(\lambda_{k_i})$ and $c_j+\lambda_{r_j}^{-1}(p_j^{n_j}A)$, and so by Corollary \ref{Corollary:IntersectionCosets} it is empty or it defines a coset of $\theta(A;0)$. This proves part (a).\\

We now prove part (b). Using an enumeration of $J$ and Lemma \ref{Lemma:intersection-kernels}, we can write $\theta(x;\bm{0})$ as 
\begin{align*}
    \theta(x;\bm{0})\,\,&:\,\,\bigwedge_{i\in I}[\,k_i\cdot x=0\,]\wedge \bigwedge_{j\in J}[\,\exists z_j(p_j^{n_j}\cdot z_j=r_j\cdot x)\,]\\
    &\equiv \exists z_1,\ldots,z_t \left((k\cdot x=0)\wedge \bigwedge_{j=1}^t[\,r_j\cdot x - p_j^{n_j}\cdot z_j=0\,]\right),
\end{align*}
where $k=\gcd(k_i:i\in I)$ by Lemma \ref{Lemma:intersection-kernels}.
We can rewrite $\theta(x;\bm{0})$ as $\exists \ov{z}(B\cdot x+C\cdot \ov{z}=\bm{0})$, using the matrix
\[M=[B\ |\ C]=\begin{bmatrix}
k&&\rvline &&0&0&\cdots&0\\
r_1&&\rvline && -p_1^{n_1} & 0 & \cdots & 0\\
r_2&&\rvline & & 0 & -p_2^{n_2} & \cdots & 0\\
\vdots &&\rvline & & \vdots & \vdots & \ddots & 0\\
r_t & &\rvline & & 0 & 0 & \cdots & -p_t^{n_t}
\end{bmatrix}\]
Given an abelian group $A$, we can consider the linear maps $M_A$ and $C_A$ corresponding to multiplication by the matrices $M$ and $C$, respectively. Let $H$ be the subgroup of $A$ defined by the formula $\theta(x;\bm{0})$, and consider the projection map $\pi:\operatorname{ker}(M_A)\to H$ given by $\pi(c;\ov{d})=c$. This is a surjective homomorphism, and its kernel is precisely $\operatorname{ker}(C_A)$.\\

Hence, we can take $\nu=\operatorname{rank}(M)$, $\mu=\operatorname{rank}(C)$, and by Lemma \ref{Lemma:KernelSNF} there are integers $c_1,\ldots,c_\nu,d_1,\ldots,d_\mu$ such that 

\begin{align*}
|H|=\dfrac{|\operatorname{ker}(M_A)|}{|\operatorname{ker}(C_A)|}=\dfrac{|A|^{|J|+1-\nu}\cdot \prod_{\ell=1}^\nu |A/c_\ell A|}{|A|^{|J|-\mu}\cdot \prod_{\ell=1}^\mu|A/d_\ell A|}=|A|^{1+\mu-\nu}\cdot \dfrac{\prod_{\ell=1}^\nu |A/c_\ell A|}{\prod_{\ell=1}^\mu|A/d_\ell A|}.
\end{align*}
Moreover, since $M$ and $B$ only differ by one column, we have $\nu\leq 1+\mu$, and so $1+\mu-\nu\in\{0,1\}$. 
\edem

\bp [Theorem 4.2.3 in \cite{AMSW}] The collection of all finite abelian groups is a m.e.c.
\label{Prop:MecAG}
\ep
\bdem By the Projection Lemma, it suffices to consider formulas of the form $\varphi(x;\ov{y})$ where $x$ is a single variable. By Lemma \ref{PruferPPF}, $\varphi(x;\ov{y})$ is equivalent to a disjunction of formulas $\tau_j(x,\overline{a})$,  where each $\tau_j(x,\overline{a})$ is a conjunction of formulas of the form ``$t(x,\ov{y})=0$'' (where $t$ is a term), or ``$p^n|t(x;\ov{y})$'' (for a prime $p$ and a positive integer $n$), and possibly with negations.\\

As usual, by increasing the number of conjunctions if necessary, we can assume that the sets defined by the formulas $\tau_i(x,\overline{a})$ are all disjoint. Consequently, we would have  $|E|=\sum_{i=1}^{k}|\tau_i(\mathcal{M},\overline{a})|$ and thus it suffices to prove the statement for sets defined by a single conjunction of basic formulas. Hence, we may suppose that $\varphi(x;\ov{y})$ has the form
\[\bigwedge_i (t_i(x;\ov{y})=0) \wedge \bigwedge_j (t_j(x;\ov{y})\neq 0)\wedge \bigwedge_r (\,p^{n_r}|t_r(x;\ov{y})\,)\wedge \bigwedge_s (\,p^{n_s}\not{|}\,\,t_s(x;\ov{y})\,).\]
We may reduce the analysis a bit further. If $\{R_i:i\in I\}$ is the collection of positive conjuncts and $\{S_j:j\in J\}$ is the collection of conjuncts appearing with negation in $\varphi(x;\ov{y})$, then by inclusion-exclusion we have:

\begin{align*}
|\varphi(x;\ov{y})|&=\left|\bigcap_{i\in I}R_i\cap \bigcap_{j\in J}\neg S_j\right|=\left|\bigcap_{i\in I}R_i\cap  \left(\bigcup_{j\in J} S_j\right)^c\right|\\
&=\left|\bigcap_{i\in I}R_i\right|-\left|\bigcap_{i\in I}R_i\cap \left(\bigcup_{j\in J} S_j\right)\right|=\left|\bigcap_{i\in I}R_i\right|-\left|\bigcup_{j\in J}\left(S_j\cap \bigcap_{i\in I}R_i\right)\right|\\
&=\left|\bigcap_{i\in I}R_i\right|-\sum_{r=1}^{|J|}(-1)^r\cdot \left(\sum_{J_0\subseteq J}\ \left|\bigcap_{j\in J_0}S_j\cap \bigcap_{i\in I}R_i\right|\ \right).
\end{align*}
Note that since the sets $I$ and $J$ are finite, to obtain the cardinality of $\varphi(x;\ov{y})$ it suffices to know the cardinalities of the intersections of \emph{positive} basic sets. Hence, it suffices to show the measurability condition for formulas of the form
\[\theta(x;\ov{y})\,\,:\,\,\bigwedge_{i\in I}(\,k_i\cdot x+t_i(\ov{y})=0\,)\wedge \bigwedge_{j\in J}(\,\,p_j^{n_j}|(r_j\cdot x+t_j(\ov{y}))\,\,)\ .\]
By Lemma \ref{Lemma:sizes-pp-definable}, for every finite abelian group $A$ the set defined by $\theta(x;\ov{a})$ is either empty, or it has cardinality
\[h_\theta:=|\theta(A;\ov{a})|=|A|^{1+\mu-\nu}\cdot \dfrac{\prod_{\ell=1}^\nu |A/c_\ell A|}{\prod_{\ell=1}^\mu |A/d_\ell A|},\]for some integers $\mu,\nu,c_1,\ldots,c_\nu,d_1,\ldots,d_\mu$ depending only on the integers $k,r_j,p_j^{n_j}$. Thus, ultimately there are only two possible sizes for $|\theta(A;\ov{a})|$ ($0$ or $h_\theta$), depending only on the formula $\theta(x;\ov{y})$, and with the definability clause given by:
\[|\theta(A;\ov{a})|=\begin{cases}h_\theta &\text{if $A\models \exists x(\theta(x;\ov{a}))$},\\ 0&\text{if $A\models \neg\exists x(\theta(x;\ov{a}))$}.\end{cases}\]
Combining this with our previous analysis, we can conclude that the definition of $\mathcal{R}$-mec holds along the class of abelian groups for every formula $\varphi(x;\ov{y})$ with a single object variable $x$, giving finitely many options depending on the cardinalities $\bigcap_{i\in I}R_i$ and $\bigcap_{i\in I}R_i\cap \bigcap_{j\in J_0}S_j$ for the possible finite subsets $J_0$ of $J$. The result then follows by the Projection Lemma (Fact \ref{Fact:Projection-Lemma}). 
\edem

\subsection{MS-measurability and pseudofiniteness of abelian groups} \label{SS:MS-measurability}

Another notion of measurability is MS-measurability in the sense of Macpherson–Steinhorn introduced in \cite{Macpherson-Steinhorn}, surveyed in \cite{Elwes-Macpherson-survey}. Roughly speaking, we say that a structure $\mathcal{M}$ is MS-measurable if there is a definable way of assigning dimensions and measures to the definable sets in $\mathcal{M}$. In work of Elwes and Macpherson (see \cite[Corollary 3.6]{Elwes-Macpherson-survey}), they showed that MS-measurable structures are supersimple of finite SU-rank. Kestner studied MS-measurability in modules, and in particular in abelian groups, providing a characterization of MS-measurable abelian groups via Szmielew groups. It is claimed in \cite[Remark 4.14]{Kestner2014} that MS-measurable abelian groups are precisely the pseudofinite abelian groups. However, in the following example we show this is false, by constructing a pseudofinite abelian group that is not MS-measurable since it fails to be superstable. First, let us recall the following Theorem due to Rogers in \cite{PhD-PatriciaRogers}. 

\bt\label{thm:Rogers} Let $A$ be an abelian group. The theory $\operatorname{Th}(A)$ is superstable if and only if the set
\[\Sigma(A)=\left\{(p,n): \text{$p$ prime, $n<\omega$, } \left|p^nA/p^{n+1}A\right|=\infty\right\}\]
is finite.
\et

With Theorem \ref{thm:Rogers} in mind, the strategy for the counterexample will be to construct a pseudofinite group $A$ for which $\Sigma(A)$ would be infinite. 

\be \label{Example:CounterexampleClaimKestner} Fix a prime $p$ and for a natural $t\geq1$, define $A_t$ as the following finite group:
\begin{align*}
    A_{t}:=\Z(p)\oplus \mathbb{Z}(p^2)\oplus\cdots\oplus \mathbb{Z}(p^t)
\end{align*}
Recall that $\text{Tf}(q,n;\Z(p^k))=p$ if $q=p$ and $n<k$, otherwise it is equal to $1$. Therefore, for all $t>n$, we obtain that $\text{Tf}(p,n;A_t)=p^{t-n}$, otherwise it is equal to $1$. Now, we consider the following ultraproduct:
\begin{align*}
    A = \prod_{t\to\mathcal{U}}A_t = \prod_{t\to\mathcal{U}}\Z(p)\oplus \mathbb{Z}(p^2)\oplus\cdots\oplus \mathbb{Z}(p^t)
\end{align*}
where we vary $t$ along $\mathcal{U}$ which is a non-principal ultrafilter on the natural numbers. Notice that $A$ is a pseudofinite abelian group, but it is not superstable because for every pair $(p,n)$ with $n\geq 0$ we have $\operatorname{Tf}(p,n;A)=\lim_{t\to\infty}\operatorname{Tf}(p,n;A_t)=\infty$. Hence, $A$ is a pseudofinite abelian group that is not MS-measurable, as this property implies supersimplicity.\\

Moreover, the ultraproduct $A=\prod_{t\to \mathcal{U}}A_t$ has the following Szmielew invariants:

\begin{table}[h!]
\centering
\begin{tabular}{|c|c|}
\hline
\textbf{Invariant} & $\prod\limits_{t\to\mathcal{U}}A_t$                                                            \\ \hline
$\operatorname{U}(q,n;-)$          & 1                                                                         \\ \hline
$\operatorname{D}(q,n;-)$          & \begin{tabular}[c]{@{}c@{}}$\infty$ if $(p=q)$\\ $1$ otherwise\end{tabular} \\ \hline
$\operatorname{Tf}(q,n;-)$         & \begin{tabular}[c]{@{}c@{}}$\infty$ if $(p=q)$\\ $1$ otherwise\end{tabular} \\ \hline
$\operatorname{Exp}(q,n;-)$        & $\infty$                                                                  \\ \hline
\end{tabular}
\caption{Szmielew invariants for the ultraproduct $A=\prod_{t\to\mathcal{U}}A_t$.}
\label{InvUltraproductoOfA_t}
\end{table}
Hence, using Table \ref{InvariantesSumaDeLosNoPseudofinitos} and Lemma \ref{BaurInvariantsForSums}, we can observe that $A\equiv (\Z(p^\infty)\oplus\Z_{(p)})^{\omega}$.
\ee
\subsection{Characterization of $N$-dimensional asymptotic classes of finite abelian groups} \label{SS:AsymptoticClasses}

Let us recall the following definition from Elwes \cite{Elwes-AsymptoticClasses}.

\bd \label{definition:N-asymptotic-classes} We say that a class $\mathcal{C}$ of finite $\Le$-structures is an \emph{$N$-dimensional asymptotic class} if for every $\Le$-formula $\varphi(x;\ov{y})$ there is a finite set of pairs $D_\varphi\subseteq (\{0,\ldots,N\}\times \mathbb{R}^{>0})\cup \{(0,0)\}$, and a finite collection of formulas $\{\psi_{(d,\mu)}(\ov{y}):(d,\mu)\in D_\varphi\}$ such that
\[\left||\varphi(M;\ov{a})|-\mu|M|^{\frac{d}{N}}\right|=o(|M|^{\frac{d}{N}}) \text{\ \ \ as $|M|\to\infty$ and $M\models \psi_{(d,\mu)}(\ov{a})$}\]
\ed

Notice that in the previous definition, if the condition holds for the pair $(0,0)$, it implies that the cardinalities $|\varphi(M;\ov{a})|$ of definable sets are uniformly finite.
\bt [Corollary 2.8 in \cite{Elwes-AsymptoticClasses}] \label{SuperstabilityAsymptoticClasses} Any ultraproduct of structures in an $N$-dimensional asymptotic class is supersimple of $\operatorname{SU}$-rank at most $N$.
\et

The concept of $N$-dimensional asymptotic classes is a natural generalization of the \emph{1-dimensional asymptotic classes}, originally defined by Macpherson and Steinhorn in \cite{Macpherson-Steinhorn}. The most important example of an asymptotic class is the class of finite fields, which is a 1-dimensional asymptotic class by a theorem of Chatzidakis, van den Dries and Macintyre (see \cite{MR1162433}). However, the most relevant example in our study would be the class of cyclic groups:

\bt [Theorem 3.14 in \cite{Macpherson-Steinhorn}] The class of finite cyclic groups is a 1-dimensional asymptotic class.
\et

The purpose of the rest of this section will be to provide conditions under which a class of finite abelian groups may or may not be an $N$-dimensional asymptotic class. We start with the following Proposition that uses non-standard cardinalities on ultraproducts of finite structures, which will be used to obtain several counterexamples.

\bp \label{Proposition:N-dimensional-mu} If $M$ is an infinite ultraproduct of structures in an $N$-dimensional asymptotic class, then for every infinite definable set $X\subseteq M^1$ there is a pair $(d,\mu)\in \{1,\ldots,N\}\times \mathbb{R}^{>0}$ such that $\operatorname{st}\left(\dfrac{|X|}{|M|^{d/N}}\right)=\mu.$
\ep
\bdem Suppose $M=\prod_{\mathcal{U}}M_n$ is an infinite ultraproduct of finite strucures in an $N$-dimensional asymptotic class $\mathcal{C}$, and let $X\subseteq M^1$ be an infinite definable set, defined by a formula $\varphi(x;\ov{a})$, with $\ov{a}=[\ov{a}_n]_\mathcal{U}\in M^{|\ov{y}|}$. By Definition \ref{definition:N-asymptotic-classes}, there is a finite set of pairs $D_\varphi$ such that
\[|\varphi(M;\ov{a})|-\mu|M|^{\frac{d}{N}}=o(|M|^{\frac{d}{N}}) \text{\ \ \ as $|M|\to\infty$ and $M\models \psi_{(d,\mu)}(\ov{a})$}\]
In particular, since $D_\varphi$ is finite, there is a single pair $(d,\mu)$ such that $M_n\models \psi_{(d,\mu)}(\ov{a}_n)$ holds for $n$ in an $\mathcal{U}$-large set $I\subseteq \mathbb{N}$. Moreover, since $X$ is infinite we know that $\mu\neq 0$, and so  we have 
\[\left(\dfrac{|\varphi(M_n;\ov{a}_n)|}{|M_n|^{d/N}}-\mu\right)\to 0,\]
as $n$ tends to infinity. By taking non-standard cardinalities and the standard part, it implies $\operatorname{st}\left(\dfrac{|X|}{|M|^{d/N}}\right)=\mu>0.$
\edem

As we have already seen, the class of finite abelian groups is a multidimensional exact class. We will explore different subclasses of finite abelian groups to see whether they correspond to $N$-dimensional asymptotic classes for some $N\in\mathbb{N}$. To start this analysis, recall that by the Fundamental Theorem of finite Abelian Groups, every finite abelian group can be written as 
\[\mathbb{Z}(p_1^{k_1})^{t_1}\oplus\cdots \oplus \mathbb{Z}(p_r^{k_r})^{t_r},\]
for some prime numbers $p_1\leq\cdots\leq p_r$, and some positive natural numbers $k_1,\ldots,k_r,t_1,\ldots,t_r$.\\

Moreover, if there is no factor of the form $\mathbb{Z}(p^k)^t$ in the decomposition, we could put $k_p=0$ and $t_p=0$, obtaining a direct sum of the form \[G\cong \bigoplus_p [\mathbb{Z}(p)^{t_{p,1}}\oplus\cdots\oplus\mathbb{Z}(p^{k_p})^{t_{p,k_p}}],\] where for all but finitely many primes we have $k_p=t_{p,j}=0$. We will call this the \emph{fundamental decomposition} of $G$.

\bd Let $p$ be a prime number, and $n\geq 1$.
\bitem
\item For a finite abelian group $G_0$ with fundamental decomposition as above, we define $p\operatorname{-width}(G_0)=\max\{t_{p,j}:j\leq k_p\}$ and $p\operatorname{-depth}(G_0)=k_p$. We also define the joint $p\operatorname{-width}\operatorname{-depth}$ value of $G_0$ as $p\operatorname{-wd}(G_0)=\min\{p\operatorname{-width}(G_0),p\operatorname{-depth}(G_0)\}$.
\item For a class of finite abelian groups, we define
\begin{align*}
p\operatorname{-width}(\mathcal{G})&:=\sup\{p\operatorname{-width}(G_0):G_0\in\mathcal{G}\} \in \mathbb{N}\cup\{\infty\},\\
p\operatorname{-depth}(\mathcal{G})&:=\sup\{p\operatorname{-depth}(G_0):G_0\in\mathcal{G}\} \in \mathbb{N}\cup\{\infty\},\\
p\operatorname{-wd}(\mathcal{G})&:=\sup\{p\operatorname{-wd}(G_0):G_0\in\mathcal{G}\} \in \mathbb{N}\cup\{\infty\}.
\end{align*}
\eitem
\ed

\bl \label{Lemma:wdinfty} Let $\mathcal{G}$ be a class of finite abelian groups. Suppose that for every $\ell\geq 1$ the class $\mathcal{G}$ contains a group $H_\ell$ with a factor of the form $\mathbb{Z}(p^{k_\ell})^{t_\ell}$ in its fundamental decomposition, with $k_\ell,t_\ell\geq \ell$. Then, there is an ultraproduct of groups in $\mathcal{G}$ that is not superstable, and $\mathcal{G}$ is not an $N$-dimensional asymptotic class for any $N\geq 1$.
\el
\bdem We can calculate the $\operatorname{Tf}$-invariant of $H_\ell$ for a pair $(p,n)$ with $n\geq 1$.

\begin{align*}
\operatorname{Tf}(p,n;H_\ell)=\prod_{i=1}^{r}\left(\operatorname{Tf}(p,n;\mathbb{Z}(p_i^{k_i})\right)^{t_i}\geq \operatorname{Tf}(p,n;\mathbb{Z}(p^{k_\ell})^{t_\ell})=\begin{cases}p^{t_\ell} &\text{if $n<k_\ell$,}\\ 1 &\text{otherwise.}\end{cases}
\end{align*}
In particular, $H_\ell\models \text{``$\operatorname{Tf}(p,n)\geq \ell$''}$ whenever $n<\ell$. Now consider the ultraproduct $A=\prod_{\mathcal{U}}H_\ell$ with respect to a non-principal ultrafilter $\mathcal{U}$. By \Los' Theorem, $A\models \text{``$\operatorname{Tf}\geq \ell$''}$ for every $n\geq 1$, and so $\operatorname{Tf}(p,n;A)=\infty$ for every pair of the form $(p,n)$, with $n\geq 1$. Therefore, by Theorem \ref{thm:Rogers}, $A$ is not superstable, and by Theorem \ref{SuperstabilityAsymptoticClasses}, this implies that $\mathcal{G}$ is not an $N$-asymptotic class for any $N\geq 1$.
\edem

Notice that in the previous Lemma the hypothesis is not simply $p-\operatorname{wd}(\mathcal{G})=\infty$, as it is necessary that both width and depth grow in the same factor. We give later an example of a class with $p-\operatorname{wd}(\mathcal{G})=\infty$ but the infinite ultraproducts are superstable of SU-rank 2. (Example \ref{Example:wd-infty-superstable}).

Observe that the use of $p-\operatorname{wd}$ in Lemma \ref{Lemma:wdinfty} gives us an asymptotic growth of both $p-\operatorname{width}$ and $p-\operatorname{depth}$ in the class. However, we are interested in softening the condition that both of them grow simultaneously, as only a controlled growth of the $p-\operatorname{width}$ in the class is enough to characterize its dimension. Therefore, we introduce the following notion of \emph{entanglement}.
\bd\label{def:Entangled} We say that a class of finite abelian groups $\mathcal{G}$ is \textbf{entangled} if given a collection of primes $\{p_\lambda\}_{\lambda<\omega}$ (possibly infinite) with $p_\lambda-\operatorname{width}(\mathcal{G})=\infty$, for all $m\in\Nat$ and every finite subset $I\subset \Nat$, there is $G_{m,I}\in\mathcal{G}$ such that for all $i\in I$, $p_{\lambda_i}-\operatorname{width}(\mathcal{G}_{m,I})\geq m$.
\ed
\brmk Note that given a class $\mathcal{G}$ where there is a single prime $p$ with $p-\operatorname{width}(\mathcal{G})=\infty$, then $\mathcal{G}$ will be entangled.
\ermk
\bej[\textbf{Entangled $1$-dimensional class with infinite $p$-width and bounded $p$-depth}] The class $\mathcal{G}=\{\mathbb{Z}(p)^t:t\geq 1\}$ has infinite $p$-width and is entangled, but it is clearly a 1-dimensional asymptotic class as all infinite ultraproducts are elementarily equivalent to $\mathbb{Z}(p)^\omega$, which is strongly minimal for every prime $p$. This example also shows that $p-\operatorname{width}(\mathcal{G})=\infty$, for a single prime $p$, is not enough to ensure that a class of finite abelian groups $\mathcal{G}$ is not an $N$-dimensional asymptotic class for some $N$.
\eej

\bl \label{lem:InfiniteWidthforInfPrimes} Let $\mathcal{G}$ be an entangled class of finite abelian groups. If $p-\operatorname{width}(\mathcal{G})=\infty$ for infinitely many primes $p$, then $\mathcal{G}$ is not an $N$-dimensional asymptotic class for any $N\geq1$. 
\el
\bdem We will construct an ultraproduct of elements in $G$ that is not superstable, let $\{p_{\lambda}\}_{\lambda<\omega}$ be the collection of infinitely many primes such that $p_\lambda-\operatorname{width}(\G)=\infty$. By the entangled property of the class, we can find a family $\{A_m\}_{m\in\Nat}\subseteq\mathcal{G}$ such that $p_{\lambda}-\operatorname{width}(A_m)\geq m$ for $\lambda\leq m$. In particular, we get the following behavior for the torsion-free invariant:
\begin{align*}
    \operatorname{Tf}(p_\lambda,0;A_{m})\geq\operatorname{Tf}\left(p_\lambda,0;\Z(p^{k_{m}})^{t_{m,i}}\right) = p^{t_{m,i}}\geq m
\end{align*}

Now, we consider the ultraproduct $A=\prod_{\mathcal{U}}A_m$ for an ultrafilter $\U$. By \Los' Theorem, we get that $\operatorname{Tf}(p_{\lambda},0;A)=\infty$ for all primes $p_{\lambda}$, and Theorem \ref{thm:Rogers} implies that $A$ is not superstable. Therefore, Theorem \ref{SuperstabilityAsymptoticClasses} implies that $\mathcal{G}$ is not an $N$-asymptotic class for any $N\geq 1$.
\edem
\bej[\textbf{Entangled classes with infinite $p$-width for infinitely many primes}] The class $\mathcal{G}=\{\mathbb{Z}(p_1)^{t_1}\oplus\cdots\oplus\mathbb{Z}(p_r)^{t_r}:r\geq 1,t_i\geq0\} $ is entangled, and for any prime $p_i$, we get $p_i-\operatorname{width}(\mathcal{G})=\infty$ and $p_i-\operatorname{depth}(\mathcal{G})=1$. Hence, by Lemma \ref{lem:InfiniteWidthforInfPrimes}, it is not an $N$-dimensional asymptotic class for any $N\geq1$. As a generalization of the previous class, for a fixed $K$ the class $\mathcal{G_K}=\{\mathbb{Z}(p_1^K)^{t_1}\oplus\cdots\oplus\mathbb{Z}(p_r^K)^{t_r}:r\geq 1,t_i\geq0\}$ is entangled and not an $N$-dimensional asymptotic class for any $N\geq 1$ by Lemma \ref{lem:InfiniteWidthforInfPrimes}. 
\eej

Note that the previous examples display situations where the collection of finite abelian groups is entangled. The following examples show the current understanding that we have about subclasses of abelian groups that may (or fail to) be $N$-dimensional asymptotic classes, and may not be entangled.

\bej[\textbf{Entangled $1$-dimensional class with infinite $p$-depth and bounded $p$-width}] For a fixed prime $p$, consider the entangled class $\mathcal{G}_p=\{\mathbb{Z}(p^k):k\geq 1\}$. If $A=\prod_{\mathcal{U}}\mathbb{Z}(p^k)$ is an infinite ultraproduct of groups in $\mathcal{G}_p$, then $\operatorname{Tf}(q,n;A)=\lim_{\mathcal{U}}\operatorname{Tf}(q,n;\mathbb{Z}(p^k))$ is either $1$ if $q\neq p$, or $p$ if $q=p$. This structure is not $\omega$-stable (Macintyre), but it is superstable. Moreover, checking the Szmielew invariants, it is elementarily equivalent to $\mathbb{Z}(p^{\infty})\oplus \mathbb{Z}_{(p)}$. 

Furthermore, $A$ is an ultraproduct of cyclic groups, and since the class of finite cyclic groups is a 1-dimensional asymptotic class we can conclude that $\operatorname{SU}(A)=1$. This example also shows that $p-\operatorname{depth}(\mathcal{G})=\infty$ is not enough to guarantee that an entangled class of finite abelian groups $\mathcal{G}$ is not an $N$-dimensional asymptotic class.
\eej

\bej[\textbf{Entangled $1$-dimensional class with $p$-width equal to $1$}] Notice that if every group in $\mathcal{G}$ has the form $\mathbb{Z}(p_1^{k_1})\oplus \cdots \oplus\mathbb{Z}(p_r^{k_r})$, with $p_1<\cdots<p_r$ distinct primes, then every group in $\mathcal{G}$ is cyclic and so it is a 1-dimensional asymptotic class, and $\mathcal{G}$ is not entangled. Hence, if $\operatorname{p-width}(\mathcal{G})=1$ for every prime $p$, the resulting ultraproduct will have SU-rank 1.
\eej
Now we will turn our interest into constructing examples of entangled classes that fail to be $N$-dimensional for $N\geq1$.
\bej [\textbf{Entangled class with infinite $p$-width and $q$-width}] Let $p,q$ be different prime numbers with $p<q$. The class $\mathcal{G}_{p,q}=\{\mathbb{Z}(p)^m\oplus\mathbb{Z}(q)^n:m,n\geq 1\}$ is \textbf{not} an $N$-dimensional asymptotic class for any $N\geq 1$. Note that it is entangled, $p-\operatorname{width}(G_{p,q})=\infty$ and $q-\operatorname{width}(G_{p,q})=\infty$, with bounded depth for any prime.\\

By Proposition \ref{Proposition:N-dimensional-mu}, it suffices to find an infinite ultraproduct $M$ of groups in $\mathcal{G}_{p,q}$ and an infinite definable set $X\subseteq M^1$ such that $\operatorname{st}\left(\dfrac{|X|}{|M|^{d/N}}\right)\neq \mu$ for any $1\leq d\leq N$ and $\mu\in\mathbb{R}^{>0}$. Let us fix $N\geq 1$ and $1\leq d\leq N$, and for every $k\geq 1$, consider the group $G_k=\mathbb{Z}(p)^{m_k}\oplus \mathbb{Z}(q)^{n_k}$ where $m_k,n_k$ are chosen to be arbitrarily large integers satisfying the following inequalities 
\[m_k\cdot \frac{d}{N}< k\text{\ \ and\ \ }n_k\cdot \left(1-\frac{d}{N}\right)>k. \text{\hspace{1cm}}(*)\]
For the case when $d=N$, the second inequality is vacuous, so we simply take a nonzero value of $n_k$ and let $m_k$ be an increasing family of indices.\\

Consider now the formula $\varphi(x):q\cdot x=0$. In the finite group $G_k=\mathbb{Z}(p)^{m_k}\oplus \mathbb{Z}(q)^{n_k}\in\mathcal{G}_{p,q}$, the formula $\varphi(x)$ defines the set $X_k:=\{0\}\oplus\mathbb{Z}(q)^{n_k}$, and $|X_k|=q^{n_k}=\dfrac{1}{p^{m_k}}|G_k|$. Hence, if $d<N$, the inequalities in $(*)$ yield the following
\begin{align*}
\lim_{k\to\infty}\dfrac{|X_k|}{|G_k|^{d/N}}=\lim_{k\to\infty} \dfrac{1}{p^{m_k}}|G_k|^{1-d/N}=\lim_{k\to\infty} \dfrac{q^{n_k\cdot (1-d/N)}}{p^{m_k\cdot d/N}}>\lim_{k\to\infty} \left(\dfrac{q}{p}\right)^k=\infty.
\end{align*}

On the other hand, if $d=N$, then we would obtain
\begin{align*}
\lim_{k\to\infty}\dfrac{|X_k|}{|G_k|^{d/N}}=\lim_{k\to\infty} \dfrac{1}{p^{m_k}}=0.
\end{align*}

In either case, this limit is not a real positive number. Hence, by taking the ultraproduct $G=\prod_{\mathcal{U}}G_k$ with respect to a non-principal ultrafilter $\mathcal{U}$ and $X=\{g\in G:q\cdot g=0\}$ we would have that $\operatorname{st}\left(\dfrac{|X|}{|G|}\right)$ is either $0$ or $\infty$. We conclude that $\mathcal{G}_{p,q}$ is not an $N$-dimensional asymptotic class, for any $N\geq 1$.
\eej

\bej[\textbf{Entangled class with equal growth of infinite $p$-width and $q$-width}]\label{ex:samepowerprimes} The non-entangled class $\mathcal{G}_2^{=}=\{\mathbb{Z}(p)^m\oplus\mathbb{Z}(q)^m:\text{$p,q$ prime, $m\geq 1$}\}$ is \textbf{not} an $N$-dimensional asymptotic class. Following the previous example, let us fix $N\geq1$, and $1\leq d\leq N$, and for every $k\geq1$, we consider the group $G_k=\mathbb{Z}(p)^{m_k}\oplus\mathbb{Z}(q)^{m_k}\in\mathcal{G}_2^{=}$ where $m_k$ satisfies $m_k> N\cdot k$, and $p<q$.\\

Consider the formula $\varphi(x):=q\cdot x=0$. In the finite group $G_k$, $\varphi(x)$ defines the set $X_k=\{0\}\oplus \Z(q)^{m_k}$, and $|X_k|=q^{m_k} = \dfrac{1}{p^{m_k}}|G_k|$. Hence, if $d=N$, we get the following growth:
\begin{align*}
    \lim\limits_{k\to\infty}\frac{|X_k|}{|G_k|}=\lim\limits_{k\to\infty}\frac{1}{p^{m_k}}=0
\end{align*}
On the other hand, for the case when $d<N$, the growth can be computed by:
\begin{align*}
    \lim\limits_{k\to\infty}\frac{|X_k|}{|G_k|^{d/N}} =\lim\limits_{k\to\infty}\frac{q^{m_k}}{(p^{m_k}q^{m_k})^{d/N}} = \lim\limits_{k\to\infty} \frac{q^{m_k(1-d/N)}}{p^{m_kd/N}}>\lim\limits_{k\to\infty}\frac{q^{k(N-d)}}{p^{kd}}\geq\lim\limits_{k\to\infty}\left(\frac{q}{p}\right)^k=\infty
\end{align*}

Hence, after taking the corresponding ultraproducts, we can conclude that $\mathcal{G}_{2}^{=}$ is not an $N$-dimensional asymptotic class, for any $N\geq1$.
\eej
We note the following nice consequence of Example \ref{ex:samepowerprimes}.
\bc If $\mathcal{G}$ is a class of finite abelian groups that is closed under direct sums, then either $\mathcal{G}\subseteq\{\mathbb{Z}(p)^t:t\geq 1\}$ (and so it is a 1-dimensional asymptotic class with ultraproducts elementarily equivalent to $\mathbb{F}_p^\omega$), or it is not an $N$-dimensional asymptotic class for any $N\geq 1$.
\ec
\bej[\textbf{Entangled class with infinite $p_i$-width for finitely many primes $p_i$}] Let $p_1<\ldots<p_M$ be different prime numbers. The entangled class $\mathcal{G}_{p_1,\ldots,p_M}=\{\mathbb{Z}(p_1)^t\oplus\cdots\oplus\mathbb{Z}(p_M)^t:t\geq 1\}$ is \textbf{not} $N$-dimensional asymptotic class for any $N\geq1$. As in the previous examples, we are going to build an ultraproduct where the standard size computation is not going to hold for a strictly positive $\mu$. So let us fix $N\geq1$ and $1\leq d\leq N$, and for every $k\geq1$, we consider the group $G_k=\Z(p_1)^{t_k}\oplus\cdots\oplus\Z(p_M)^{t_k}\in\mathcal{G}_{p_1,\ldots,p_M}$, where $t_k$ satisfies $t_k>N\cdot k$ if $d=N$. Otherwise, $t_k$ is given by the constraint:
\begin{align*}
    t_k>\max\left\{\dfrac{(M-1)\cdot k}{(1-d/N)},k\cdot N/d\right\}
\end{align*}

Similarly, we will consider the formula $\varphi_M(x):=p_M\cdot x=0$. In our finite group $G_k$, $\varphi_M$ defines the set $X_k = \{0\}\oplus\cdots\oplus\Z(p_M)^{t_k}$ which has size $p_M^{t_k}=\dfrac{|G_k|}{(p_1\cdots p_{M-1})^{t_k}}$. Then we can compute our size estimates as follows, let us begin with the case $d=N$:
\begin{align*}
    \lim\limits_{k\to\infty}\frac{|X_k|}{|G_k|} = \lim\limits_{k\to\infty}\frac{1}{(p_1\cdots p_{M-1})^{t_k}}\to0
\end{align*}

On the other hand, for the case that $d<N$, we obtain the following bound:
\begin{align*}
    \lim\limits_{k\to\infty}\frac{|X_k|}{|G_k|^{d/N}} = \lim\limits_{k\to\infty}\frac{p_M^{t_k}}{(p_1\cdots p_M)^{t_k d/N}}=\lim\limits_{k\to\infty}\frac{p_M^{t_k(1-d/N)}}{(p_1\cdots p_{M-1})^{t_k d/N}} > \lim\limits_{k\to\infty}\left(\frac{p_M^{(M-1)}}{p_1\cdots p_{M-1}}\right)^{k}\to\infty
\end{align*}
\eej
Now, we give a positive example where we attain asymptotic dimensionality for a non-entangled class of finite abelian groups.
\bej[\textbf{Entangled $2$-dimensional class with infinite $p$-width and $p$-depth$=2$}]\label{Ex:G_p,2} For a prime $p$, note that the class $\mathcal{G}_{p,2}=\{\mathbb{Z}(p^2)^t:t\geq 1\}$ is entangled and satisfies that all its ultraproducts are equivalent to a group of the form $\Z(p^2)^{\lambda}$ where $\lambda\leq\omega$. By Macintyre's Theorem (see Theorem \ref{TeoremaMacintyre2.0}), every ultraproduct of this family is $\omega$-stable as it has bounded exponent. Furthermore, we will show that $\mathcal{G}_{p,2}$ is $2$-dimensional asymptotic.\\

First, it is easy to check that it is not $1$-dimensional asymptotic since the formula $\varphi_p(x):px=0$, that defines the $p$-torsion of each group, satisfies the following growth:
\begin{align*}
    \left|\varphi_p\left(\Z(p^2)^t\right)\right| = p^t = \sqrt{|\Z(p^2)^t|}
\end{align*}
Hence, we can consider $G_k:=\Z(p^2)^k$ and the definable set $X_k$ by the formula $\varphi_p$ on the corresponding $G_k$. Then:
\begin{align*}
    \lim\limits_{k\to\infty}\frac{|X_k|}{|G_k|}=\lim\limits_{k\to\infty}\frac{1}{\sqrt{|G_k|}}=\lim\limits_{k\to\infty}\frac{1}{p^k}=0
\end{align*}
and $|X_k|\to\infty$ as $k$ increases. Thus, $\mathcal{G}_{p,2}$ is not $1$-dimensional asymptotic. In order to show that the class is $2$-dimensional asymptotic, we will use the results from Proposition \ref{Prop:MecAG}. In particular, we will reduce our computation to only positive primitive formulas and use the result of the proposition to conclude that given a positive primitive formula $\theta(x,\overline{y})$, and Lemma \ref{Lemma:sizes-pp-definable} yields the following quantity for the corresponding $\mu,\nu$ and $c_\ell,d_\ell$ that only depend on the structure of $\theta$:
\begin{align*}
    |\theta(\Z(p^2)^t;\ov{a})|=|\Z(p^2)^t|^{1+\mu-\nu}\cdot \dfrac{\prod_{\ell=1}^\nu |\Z(p^2)^t/c_\ell \Z(p^2)^t|}{\prod_{\ell=1}^\mu|\Z(p^2)^t/d_\ell \Z(p^2)^t|}
\end{align*}
However, since $\theta$ is a positive primitive formula, we can simplify our size calculation to the following:
\begin{align*}
    |\Z(p^2)^t/c_\ell \Z(p^2)^t| = |\Z(p^2)/c_\ell \Z(p^2)|^t  = \begin{cases}
        p^{2t}, & \gcd(c_\ell,p^2)=p^2\\
        p^t, & \gcd(c_\ell,p^2)=p\\
        1, & \gcd(c_\ell,p^2)=1
    \end{cases}
\end{align*}
Then, we can rewrite $\nu = \nu_0+\nu_1+\nu_2$ where $\nu_i=\left|\{\ell:(c_\ell,p^2)=p^{2-i}\}\right|$. Similarly, we write $\mu=\mu_0+\mu_1+\mu_2$ for the collection of constants $d_\ell$. Using the previous constants, we can describe the size of the definable sets by:
\begin{align*}
    |\theta(\Z(p^2)^t;\ov{a})|&=|\Z(p^2)^t|^{1+\mu-\nu}\cdot \dfrac{\prod_{\ell=1}^\nu |\Z(p^2)^t/c_\ell \Z(p^2)^t|}{\prod_{\ell=1}^\mu|\Z(p^2)^t/d_\ell \Z(p^2)^t|}\\
    &=(p^{2t})^{1+\mu-\nu}\cdot \dfrac{(p^{2t})^{\nu_0}\cdot(p^{t})^{\nu_1}\cdot(1)^{\nu_2}}{(p^{2t})^{\mu_0}\cdot(p^{t})^{\mu_1}\cdot(1)^{\mu_2}}\\
    &=(p^{2t})^{1+\mu-\nu}\cdot \dfrac{(p^{2t})^{\nu_0}\cdot(p^{t})^{\nu_1}}{(p^{2t})^{\mu_0}\cdot(p^{t})^{\mu_1}}\\
    &=(p^{2t})^{1+\mu_2-\nu_2}\cdot (p^{2t})^{(\mu_1-\nu_1)/2}=|\Z(p^2)^t|^{e_{\theta}}.
\end{align*}
Where $e_{\theta} = 1+\mu_2-\nu_2 + \frac{\mu_1-\nu_1}{2}$, and it only depends on $\theta$. Hence, we reduce the size computation to one of the following three asymptotic behaviors:
\begin{align*}
    |\theta(\Z(p^2)^t;\ov{a})|=\begin{cases}
        0& ,\Z(p^2)^t\not\models\exists x(\theta(x,\overline{a}))\\
        |\Z(p^2)^t|^{e_{\theta}} & ,\Z(p^2)^t\models\exists x(\theta(x,\overline{a}))
    \end{cases}
\end{align*}
Therefore, we can pick $(\mu,d)$ depending on the value of $e_\theta$ to characterize its asymptotic behavior. So, $\mathcal{G}_{p,2}$ is a $2$-dimensional asymptotic class.
\eej
The following example generalizes the asymptotic class presented before.
\bej [\textbf{Entangled $N$-dimensional class with infinite $p$-width and $p$-depth$=N$}] For a prime $p$ and $M\geq1$, note that the class $\mathcal{G}_{p,k}=\left\{\mathbb{Z}(p^N)^t:t\geq 1\right\}$ is entangled, $N$-dimensional asymptotic and every ultraproduct of this family is $\omega$-stable. 
\eej

For the case of entangled classes having a single prime with infinite $p-\operatorname{width}$, we do not claim that there is a bound for the dimension depending on its depth as we will have the following examples that display a different behavior.

\bej [\textbf{Entangled $2$-dimensional class with infinite $p$-width and $p$-depth}]\label{Example:wd-infty-superstable} Consider the entangled class $\mathcal{G}=\{B_t:t\geq 1\}$, where $B_t=\mathbb{Z}(p)^t\oplus \mathbb{Z}(p^t)$. Notice that $p-\operatorname{wd}(\mathcal{G})=\infty$, and we can compute the torsion free invariant for any element in the class as follows by letting $q$ be a prime and $n\geq0$
\begin{align*}
    \operatorname{Tf}(q,n;B_t)=\operatorname{Tf}(q,n;\Z(p)^t)\cdot\operatorname{Tf}(q,n;\Z(p^t))=\begin{cases}
        p^{t+1},&n=0,q=p\\
        p,&0<n<t,q=p\\
        1,&\text{otherwise}
    \end{cases}
\end{align*}
Let $A$ be an infinite ultraproduct of groups in $\mathcal{G}$, by the previous computation, we get that the torsion free invariants of $A$ are given by:
\begin{align*}
    \operatorname{Tf}(q,n;A)=\begin{cases}
        \infty,&n=0,q=p\\
        p,&n\neq0,q=p\\
        1,&\text{otherwise}
    \end{cases}
\end{align*}
By Theorem \ref{thm:Rogers}, $A$ is superstable. When we study the dimensionality of the class $\mathcal{G}$, we can consider the definable set by formula $\varphi_p(B_t) = \Z(p)^t\oplus p^{t-1}\Z(p^t)$, thus $|\varphi_p(B_t)|=p^t\cdot p=p\sqrt{|B_t|}$, which seems to work asymptotically for $2$-dimensional. If you consider $\varphi_{p^\alpha}$ for $\alpha$ arbitrary, $|\varphi_{p^{\alpha}}(B_t)|=p^\alpha\sqrt{|B_t|}$. Therefore, by a similar procedure as in Example \ref{Ex:G_p,2}, one can show that $\mathcal{G}$ is a $2$-dimensional asymptotic class.
\eej

\bej[\textbf{Entangled $N+1$-dimensional class with infinite $p$-width and $p$-depth$=N$}]\label{Ex:G_Nnondepth} For a fixed $N>1$, consider the entangled class $\mathcal{G}_{N}=\{H_{N,t}:t\geq1\}$, where $H_{N,t} = \Z(p)^t\oplus\Z(p^N)^t$. For this class $p-\operatorname{width}(\mathcal{G}_N)=\infty$ but $p-\operatorname{depth}(\mathcal{G}_N)=N$. First, we can compute the torsion free invariants for elements in the class as follows:
\begin{align*}
    \operatorname{Tf}(q,n;H_{N,t})=\operatorname{Tf}(q,n;\Z(p)^t)\cdot\operatorname{Tf}(q,n;\Z(p^N)^t)=\begin{cases}
        p^{2t},&n=0,q=p\\
        p^{t},&0<n<N,q=p\\
        1,&\text{otherwise}
    \end{cases}
\end{align*}
where $q$ is a prime and $n\geq0$. Hence, for any infinite ultraproduct $A$ of groups in $\mathcal{G}_N$, we get the following torsion free invariants:
\begin{align*}
    \operatorname{Tf}(q,n;A)=\begin{cases}
        \infty,&n<N,q=p\\
        1,&\text{otherwise}
    \end{cases}
\end{align*}
As before, by Theorem \ref{thm:Rogers}, $A$ is superstable. Moreover, we will compute the growth of the definable set by the formula $\varphi_\alpha(x):=p^\alpha\cdot x=0$ as $t$ goes to infinity with respect to the size of $H_{N,t}$, where $\alpha\in\Nat$. Observe that if $\alpha\geq N$, $\varphi_\alpha$ defines all the group $H_{N,t}$ for any $t$, and the quotient of the sizes will be one. So, we will assume that $1\leq\alpha<N$ as the $\alpha=0$ case defines a single point. Now, we will compute the size of the definable set as follows:
\begin{align*}
    |\varphi_\alpha(H_{N,t})|=|\Z(p)^t\oplus p^{N-\alpha}\Z(p^N)^t| = p^t\cdot p^{\alpha t} = p^{(\alpha+1)t}.
\end{align*}
As $|H_{N,t}|=p^{(N+1)t}$, we compute the asymptotic growth of the ratio by the following limit:
\begin{align*}
    \lim\limits_{t\to\infty}\frac{|\varphi_\alpha(H_{N,t})|}{|H_{N,t}|} = \lim\limits_{t\to\infty}\frac{p^{(\alpha+1)t}}{p^{(N+1)t}} = \lim\limits_{t\to\infty}\frac{1}{p^{(N-\alpha)t}}=0.
\end{align*}
Since $|\varphi_\alpha(H_{N,t})|=p^{(\alpha+1)t}\to\infty$ as $t$ grows, we can conclude that $\mathcal{G}_{N}$ is not a $1$-dimensional asymptotic class. Nevertheless, we want to observe the following for the previous example. If we investigate if the class is $M$-dimensional for some $M\geq1$, we would start by finding the correct dimension for all the sets defined by $\varphi_\alpha$. So, assuming that there is a finite positive dimension $\gamma=\frac{d}{M}$ for some $d\leq M$ and $M$ finite, we would get the following:
\begin{align*}
    \lim\limits_{t\to\infty}\frac{|\varphi_\alpha(H_{N,t})|}{|H_{N,t}|^{\gamma}} = \lim\limits_{t\to\infty}\frac{p^{(\alpha+1)t}}{p^{(N+1)t\gamma}} = \lim\limits_{t\to\infty} p^{t\left((\alpha+1)-\gamma(N+1)\right)},
\end{align*}
which is finite and non-zero, only when the exponent vanishes. Hence, we get the following description for $\gamma$:
\begin{align*}
    0=(\alpha+1)-\gamma(N+1)\Longleftrightarrow\gamma=\frac{d}{M}=\frac{\alpha+1}{N+1}
\end{align*}
In particular, this implies that $M\equiv N+1\mod (\alpha+1)$, for all $0<\alpha<N$. Thus, $M=N+1$ and we can simply take $d=\alpha+1\leq N$ to obtain the correct asymptotic growth. By a similar approach as in Example \ref{Ex:G_p,2}, one can show that the class $\mathcal{G}_N$ is a $N+1$-dimensional asymptotic class.
\eej

Examples \ref{Example:wd-infty-superstable} and \ref{Ex:G_Nnondepth} are two examples of asymptotic classes of finite abelian groups where the asymptotic dimension is not equal to the depth with respect to the corresponding prime with infinite width. After observing the examples of entangled classes, we observe that the most relevant behavior that causes the asymptotic-ness to fail comes from the construction of \emph{entangled} sets that witness the growth of certain $p$-groups. However, if we restrict to classes that do not have this property of entanglement, we conjecture the following dichotomy:

\bcon\label{Conj:NDimAsympt} Let $\mathcal{G}$ be a class of finite abelian groups such that for all primes $p$, $p-\operatorname{width}(\mathcal{G})$ is finite. Then either $\mathcal{G}$ is a $1$-dimensional asymptotic class, or $\mathcal{G}$ is not $N$-dimensional asymptotic for any $N\geq1$.
\econ

\section{Appendix: SU-rank in pseudofinite abelian groups}

In this appendix, we characterize forking and SU-rank in pseudofinite abelian groups, in the language $\Le=\{+,0,(\lambda_m)_{m\in\mathbb{Z}}\}$. To start, it is easy to show that if $(A,+)$ is an abelian group then $\operatorname{Th}(A,+)$ is stable. For instance, one can show that the basic formulas $t(x;\ov{a})=0$ and $p^n|t(x;\ov{a})$ do not have the order property, and so quantifier elimination would imply that no formula in one variable would have the order property. So, the theory of $(A,+)$ is stable.\\

The original proof of Rogers' Theorem (see Theorem \ref{thm:Rogers}), which characterizes whether an abelian group $A$ is superstable, relies on the analysis of the number of types for Szmielew groups, showing that whenever $\Sigma$ is finite then the number of types over a set $B$ of cardinality $\kappa$ satisfies $S_1(B)=\kappa+2^{\aleph_0}$, whereas if $\Sigma(A)$ is infinite there are subsets $B$ of size $\kappa$ such that $S_1(B)=\kappa^{\aleph_0}$. This is a well-known dividing line on the stability function, due to Shelah (cf. \cite{Keisler-SixClasses}). \\

We will prove the same theorem for pseudofinite abelian groups via a finer analysis of forking and SU-rank for arbitrary abelian groups, which in turn will give simultaneously a characterization of forking and SU-rank. As $\operatorname{Th}(A,+)$ is stable and dividing is equivalent to forking for stable theories, we will only need to understand the dividing formulas in one variable. 

\bl \label{Lemma:decreasing-chain} Let $A$ be an infinite abelian group, and suppose $\Sigma(A)=\{(p_i,n_i):1\leq i<\kappa\}$ is an enumeration in lexicographical order of $\Sigma(A)$, with $\kappa\in\mathbb{N}\cup\{\omega\}$. Define recursively the sequence of integers $\ell_0=1$, and $\ell_i=\operatorname{lcm}(p_1^{n_1+1},\ldots,p_i^{n_{i}+1})$, and consider the sequence of groups defined as $G_i=\ell_i\cdot A$.

Then, $\{G_i:0\leq i<\kappa\}$ is a decreasing sequence of definable groups satisfying $[G_i:G_{i+1}]=\infty$ for all $0\leq i<\kappa$. 
\el
\bdem Notice that $\ell_i$ divides $\ell_{i+1}$, so $G_i\geq G_{i+1}$ for every $0\leq i<\kappa$. It is also clear that each $G_i$ is a definable subgroup of $A$. Furthermore, $G_0=A$ and $G_1=p_1^{n_1+1}A$, and since $(p_i,n_i)\in \Sigma(A)$ we have $|G_0/G_1|=|A/p_1^{n_1+1}A|\geq |p_1^{n_1}A/p_1^{n_1+1}A|=\infty$.\\

Let us show now that $|G_i/G_{i+1}|=\infty$ for every $1\leq i<\kappa$. We prove this in two cases. First, if $p_i\neq p_{i+1}$, hence $\ell_{i+1}=\operatorname{lcm}(\ell_i,p_{i+1}^{n_{i+1}})=p_{i+1}^{n_{i+1}}\cdot \ell_i$ and $\gcd(\ell_i,p_{i+1}^{n_{i+1}})=1$. Since $(p_{i+1},n_{i+1})\in\Sigma(A)$, we have $|p_{i+1}^{n_{i+1}}A/p_{i+1}^{n_{i+1}+1}A|=\infty$, and since $\gcd(\ell_i,p_{i+1})=1$, we have by Lemma \ref{Lemma:intersection-cosets} that the coset $G_i=\ell_iA$ intersects every coset of $p_{i+1}^{n_{i+1}+1}A$ contained in $p_{i+1}^{n_{i+1}}A$, and in particular, \[|G_i/p_{i+1}^{n_{i+1}+1}G_i|=|G_i/p_{i+1}^{n_{i+1}+1}\ell_i A|=|G_i/\ell_{i+1}A|=|G_i/G_{i+1}|=\infty.\]
On the other hand, if $p_i=p_{i+1}=p$, then we have $\ell_{i+1}=p^{n_{i+1}+1}\cdot r$ and $\ell_i=p^{n_i+1}\cdot r$ for some integer $r$ with $\gcd(p,r)=1$. We also have $|p^{n_{i+1}}A/p^{n_{i+1}+1}A|=\infty$, and since $(p,n_i)<_{\operatorname{lex}}(p,n_{i+1})$, we must have $n_{i+1}\geq n_i+1$. Hence, by Lemma \ref{Lemma:intersection-cosets} the coset $rA$ intersects every coset of $p^{n_{i+1}+1}A$ inside $p^{n_{i+1}}A$, and using Lemma \ref{Lemma:intersections-lcm} we obtain
\begin{align*}
|G_i/G_{i+1}|=|\ell_i\cdot A/\ell_{i+1}\cdot A|=|p^{n_i+1}r\cdot A/p^{n_{i+1}+1}r\cdot A|\geq |p^{n_i+1}r\cdot A/p^{n_{i+1}+1}r\cdot A|=\infty.
\end{align*}
Hence, $\{G_i:0\leq i<\kappa\}$ is a decreasing sequence of definable groups, with $[G_i:G_{i+1}]=\infty$ for each $0\leq i<\kappa$.
\edem

\bt \label{Theorem:TFinfty-not-superstable} Let $A$ be an infinite abelian group. If $\Sigma(A)$ is infinite then $\operatorname{Th}(A)$ is not superstable.
\et
\bdem We may assume that $A$ is saturated, and let $\alpha$ be a generic element of $A$. Since $\Sigma(A)$ is infinite, by Lemma \ref{Lemma:decreasing-chain} there exists a decreasing sequence of definable subgroups $\{G_i:i<\omega\}$ with $G_0=A$ and $[G_i:G_{i+1}]=\infty$ for every $i<\omega$. By compactness, we can choose elements $\{b_i:i<\omega\}$ such that $\alpha-b_i\in G_i\setminus G_{i+1}$ for every $i<\omega$, and by Erd\H{o}s-Rado we may even assume that $\{b_i:i<\omega\}$ forms an indiscernible sequence. To contradict superstability, we will show that there is no finite subset $B_0\subseteq \{b_i:i<\omega\}$ such that $\alpha\ind_{B_0}\{b_i:i<\omega\}$.\\

If $B_0$ is a finite subset of $\{b_i:i<\omega\}$, there is some $n<\omega$ such that $B_0\subseteq \{b_0,\ldots,b_n\}$. So, the formula $(x-b_{n+1})\in G_{n+1}$ belongs to the type $\operatorname{tp}(\alpha/B)$, and it forks over $B_0$ because there are elements $\{b_{n+1}^j:j<\omega\}$ in different $G_{n+1}$-cosets with the same type of $b_{n+1}$ over $B$, and so the set $\{(x-b_{n+1}^j)\in G_{n+1}:j<\omega\}$ would be 2-inconsistent.\\

Hence, there is no finite subset $B_0$ of $B=\{b_i:i<\omega\}$ such that $\alpha\ind_{B_0}B$, which implies that $\operatorname{Th}(A)$ is not superstable.
\edem

\bd Let $A$ be an infinite abelian group with $\Sigma(A)=\{(p_1,n_1),\ldots,(p_k,n_k)\}$ finite. We put $t=\operatorname{lcm}(p_1^{n_1+1},\ldots,p_k^{n_k+1})$, and define the invariant $\rho(A)$ given by \[\rho(A)=\begin{cases}|\Sigma(A)|+1&\text{if $|t\cdot A|=\infty$}\\ |\Sigma(A)|&\text{if $|t\cdot A|<\infty$}.\end{cases}\]
\ed

\bt \label{Theorem:CharacterizationSUrank} Let $A$ be an infinite abelian group, and suppose $\Sigma(A)$ is finite. Then, $\operatorname{Th}(A)$ is superstable with $\operatorname{SU}(A)\geq \rho(A)$. Moreover, if $A$ is pseudofinite, we have the equality $\operatorname{SU}(A)=\rho(A)$.
\et

We will need first a couple of lemmas regarding forking in arbitrary abelian groups. The following is an easy exercise in forking for modules, and corresponds precisely to Proposition 5.14 in \cite{Prest-MTandModules}.

\bh \label{Lemma:dividing-pp-sets} Let $A$ be an infinite abelian group and suppose $\varphi(x)$ is a formula defined as $\theta(\ov{x};\ov{a})\wedge \bigwedge_{j=1}^n\neg\psi_j(\ov{x};\ov{b}_j)$, where all the formulas $\theta(\ov{x};\ov{a})$ and $\psi_j(\ov{x};\ov{b}_j)$ are positive primitive formulas. Then, $\operatorname{SU}(\varphi(x))=\operatorname{SU}(\theta(x;\ov{a}))$. 
\eh

\bl \label{Lemma:dividing-chain-groups} Let $A$ be an infinite superstable abelian group, and let $p_0\subseteq p_1\subseteq \cdots \subseteq p_d$ be a dividing chain of complete types (over different sets of parameters) of length $d$. Then, there is a sequence of $\emptyset$-definable groups $H_0\geq H_1\geq \cdots \geq H_d$ with $[H_i:H_{i+1}]=\infty$ such that for every $i$ we have $\operatorname{SU}(p_i)=\operatorname{SU}(H_i)$.
\el
\bdem For $i\leq d$, let $B_i$ be the collection of parameters for the type $p_i$, so that $B_0\subseteq B_1\subseteq \cdots \subseteq B_d$. Hence, for every $i=1,\ldots,d$, there is a formula $\varphi(x;\ov{b}_i)\in p_i$ which divides over $B_{i-1}$, and so that $\operatorname{SU}(p_i)=\operatorname{SU}(\varphi(x;\ov{b}_i))$. Since we are considering complete types, we may assume by quantifier elimination that each of the formulas is a conjunction of positive primitive formulas and their negations. Hence, every formula $\varphi(x;\ov{b}_i)$ is equivalent to a formula of the form $\theta_i(x;\ov{b}_i)\wedge \bigwedge_{j=1}^n \neg\psi_{i,j}(x;\ov{b}_i)$ where $\theta(x;\ov{y})$ and $\psi_j(x;\ov{y})$ are all positive primitive. \\

Moreover, by taking conjunctions with previous formulas, we may even assume that both $\{\varphi(A;\ov{b}_i):i\leq d\}$ and $\{\theta_i(x;\ov{b}_i):i\leq d\}$ form decreasing chains of definable groups. So, by Fact \ref{Lemma:dividing-pp-sets}, we have that $\operatorname{SU}(p_i)=\operatorname{SU}(\varphi(x;\ov{b}_i))=\operatorname{SU}(\theta_i(x;\ov{b}_i))$. Moreover, the formula $\theta_i(x;\ov{b}_i)$ defines in $A$ a coset of the form $c_i+H_i$, where $H_i=\theta_i(A;\bm{0})$ is a $\emptyset$-definable group, and of course we have $\operatorname{SU}(\theta_i(x;\ov{b}_i))=\operatorname{SU}(c_i+H_i)=\operatorname{SU}(H_i)$ because the SU-rank is preserved under definable bijections.\\

Therefore, we have obtained a sequence of definable groups $H_0\geq H_1\geq \cdots \geq H_d$ with $\operatorname{SU}(p_i)=\operatorname{SU}(\theta_i(\ov{x};\ov{b}_i))=\operatorname{SU}(H_i)$ for every $i\leq d$, and since $p_{i+1}$ is a dividing extension of $p_i$, we must have $SU(H_{i+1})=\operatorname{SU}(p_{i+1})<\operatorname{SU}(p_i)=\operatorname{SU}(H_i)$, which implies that $[H_i:H_{i+1}]=\infty$.
\edem

Now we can proceed to the proof of the characterization of the SU-rank for pseudofinite abelian groups.

\bdem[Proof of Theorem \ref{Theorem:CharacterizationSUrank}]

We show first the inequality $\operatorname{SU}(A)\geq \rho(A)$ for arbitrary abelian groups. Let $\alpha$ be a generic element of $A$, and consider the finite sequence of definable groups $G_0=A\geq G_1\geq \cdots \geq G_k$ given by Lemma \ref{Lemma:decreasing-chain}, and so we have $[G_i:G_{i+1}]=\infty$ for $0\leq i<k$ and \[G_k=\operatorname{lcm}(p_1^{n_1+1},\ldots,p_k^{n_k+1})\cdot A=t\cdot A.\] Consider the sequence of partial types given by $p_0=\{x=x\}$ and $p_{i+1}=p_i\cup \{(x-b_i)\in G_i\}$ for $0\leq i\leq k$. Using an argument similar to the proof of Theorem \ref{Theorem:TFinfty-not-superstable} we have that $p_{i+1}$ is an extension of $p_i$ that divides over $B_i=\{b_0,\ldots,b_i\}$, and so $p_0\subseteq p_1\subseteq \cdots \subseteq p_k$ is a dividing chain of length $k$. Now, if $t\cdot A$ is finite, then $p_k\supseteq \{(x-b_k)\in t\cdot A\}$ is already algebraic. Otherwise, we may extend $p_k$ further by taking the type $p_{k+1}=\{x=\alpha\}$. This shows that $\operatorname{SU}(A)$ is at least $k=|\Sigma(A)|$ if $t\cdot A$ is finite, or at least $k+1$ if $t\cdot A$ is infinite, and so $\operatorname{SU}(A)\geq \rho(A)$\\

Now we prove the other inequality for pseudofinite abelian groups. Suppose $\operatorname{SU}(A)\geq \rho(A)+1$. By Lemma \ref{Lemma:dividing-chain-groups}, this implies that there is a sequence of $\emptyset$-definable groups $H_0>H_1>\cdots>H_{\rho(A)+1}$ with $[H_i:H_{i+1}]=\infty$ for every $i\leq \rho(A)$. These groups are defined by positive primitive formulas, so each of the groups can be written as an intersection of the form
\[\operatorname{ker}(\lambda_m)\cap \bigcap_{j=1}^n \{x\in A: r\cdot x\in p^nA\},\]
for integers $m,r$, $n\geq 1$ and $p$ prime. By Lemmas \ref{Lemma-CountingKernel} and \ref{Lemma-CountingDivisibility}, the cardinalities of each of the groups $H_i$ is a finite product of non-standard cardinalities of the form $|p^nA/p^{n+1}A|$. So, whenever $[H_i:H_{i+1}]=\infty$, there must be at least a pair $(p_i,n_i)\in \Sigma(A)$ witnessing $|H_i/H_{i+1}|\geq |p_i^{n_i}A/p_i^{n_i+1}A|=\infty$. Moreover, since $H_i\geq H_{i+1}$, none of the pairs $(p_i,n_i)$ can appear twice in the sequence.\\

This implies that the sequence $H_0>H_1>\cdots>H_{|\Sigma(A)|}$ has already use all the pairs $(p,n)\in \Sigma(A)$. We have two cases: if $|\Sigma(A)|=\rho(A)$, it implies that $H_{\rho(A)}\subseteq tA$ and so it is finite, contradicting $[H_{\rho(A)}:H_{\rho(A)+1}]=\infty$. On the other hand, if $|\Sigma(A)|=\rho(A)-1$, it implies that $H_{\rho(A)-1}\subseteq tA$ and both the quotient  $H_{\rho(A)-1}/H_{\rho(A)}$ and the subgroup $H_{\rho(A)}$ are infinite, but as we have seen this can only occur via a pair $(p,n)\in \Sigma(A)$ and we have already exhausted the list. So, in either case we obtain a contradiction, and we conclude that $SU(A)\leq \rho(A)$, as desired
\edem

\bej Consider the group $A=\mathbb{Z}(p^\infty)^\omega$ which is not pseudofinite. Since $A$ is divisible, we have $p^nA=A$ for every prime power $p^n$, so $\operatorname{Tf}(p,n;A)=1$ for every pair $(p,n)$, and $\Sigma(A)=\emptyset$. This implies by Theorem \ref{thm:Rogers} that $\operatorname{Th}(A)$ is superstable, but it does not have finite SU-rank because for every $N\geq 1$ we can consider the decreasing sequence of subgroups given by $H_i=\operatorname{ker}(\lambda_{p^i})$, obtaining $A\geq H_N\geq H_{N-1}\geq \cdots \geq H_1\geq 0$, and $[H_i:H_{i-1}]=|(\mathbb{Z}/p\mathbb{Z})^\omega|=\infty$. This shows that $\operatorname{SU}(A)\geq \omega$. 
\eej
\bc Let $A$ be an infinite pseudofinite abelian group. Then, $\operatorname{Th}(A)$ is either superstable of finite SU-rank, or $\operatorname{Th}(A)$ is strictly stable.
\ec
\bdem By Theorem \ref{Theorem:TFinfty-not-superstable}, if $\Sigma(A)$ is infinite then $\operatorname{Th}(A)$ is not superstable, hence it is strictly stable. On the other hand, by Theorem \ref{Theorem:CharacterizationSUrank}, if $\Sigma(A)$ is finite then the SU-rank of $A$ is bounded by $|\Sigma(A)|+1$, and so it is finite.
\edem

\bibliographystyle{plain}
\bibliography{pseudo_groups.bib}

\end{document}